\documentclass[11pt,a4paper]{article}

\newtheorem{problem}{Problem}
\newtheorem{remark}{Remark}

\usepackage{amssymb}
\usepackage{amsfonts}
\expandafter\let\csname equation*\endcsname\relax
\expandafter\let\csname endequation*\endcsname\relax
\usepackage{amsmath}

\usepackage{graphicx}

\begin{document}
\date{}

\title{Order Reduction of the Radiative Heat Transfer Model
for the Simulation of Plasma Arcs}
\author{L. Fagiano,
 R. Gati
\thanks{L. Fagiano and R. Gati are with ABB Switzerland Ltd., Corporate Research, Segelhofstrasse 1K, Baden-Daettwil, CH-5405. Corresponding author: L. Fagiano (lorenzo.fagiano@ch.abb.com)}}
\maketitle

\begin{abstract}
An approach to derive low-complexity models describing thermal
radiation for the sake of simulating the behavior of electric arcs
in switchgear systems is presented. The idea is to approximate the
(high dimensional) full-order equations, modeling the propagation of
the radiated intensity in space, with a model of much lower
dimension, whose parameters are identified by means of nonlinear
system identification techniques. The low-order model preserves the
main structural aspects of the full-order one, and its parameters
can be straightforwardly used in arc simulation tools based on
computational fluid dynamics. In particular, the model parameters
can be used together with the common approaches to resolve radiation
in magnetohydrodynamic simulations, including the discrete-ordinate
method, the P-N methods and photohydrodynamics. The proposed order
reduction approach is able to systematically compute the
partitioning of the electromagnetic spectrum in frequency bands, and
the related absorption coefficients, that yield the best matching
with respect to the finely resolved absorption spectrum of the
considered gaseous medium. It is shown how the problem's structure
can be exploited to improve the computational efficiency when
solving the resulting nonlinear optimization problem. In addition to
the order reduction approach and the related computational aspects,
an analysis by means of Laplace transform is presented, providing a
justification to the use of very low orders in the reduction
procedure as compared with the full-order model. Finally,
comparisons between the full-order model and the reduced-order ones
are presented.
\end{abstract}

%
\vspace{2pc} \noindent{\it Keywords}: Arc simulations, Radiative
heat transfer, Model order reduction, Nonlinear estimation,
Nonlinear model identification
%
%
%
%

\section{Introduction}\label{S:intro}

The switching performance of circuit breakers depends strongly on
the behavior of the electric arc that originates when the contacts
are opened in presence of relatively large electric current values
\cite{SSKPJ14,BHEL05}. In turn, the arc dynamics are influenced by
multiple interacting physical phenomena which, together with the
short timescale of the arcing event and the large values of
temperature and pressure, increase the complexity and difficulty of
understanding, carrying out experiments, and deriving numerical
models of the switching behavior. Computational fluid dynamic (CFD)
approaches are being used in both public and private research
efforts to simulate the time evolution of the plasma that carries
the current during the interruption process, see e.g.
\cite{ClTh95,BaSV08,AnLu09,ItTZ11}. The CFD simulations are often
coupled with solvers for the electro-magnetic (EM) phenomena,
resulting in sophisticated multi-physics simulation tools (see
\cite{BaSV08,ItTZ11}) that allow one to gather an insight of what is
actually happening during the current interruption process - aspects
that are very difficult to quantify with direct measurements for the
above-mentioned reasons. Such simulation tools provide a significant
added value to explain the results of experimental tests and to
support the development of switchgear devices, however they also
bring forth an important issue in addition to the inherent
difficulty of plasma modeling: the need to find a good balance
between the accuracy of the employed physical models and their
computational complexity. The modeling of radiative heat transfer
during the arcing process is an illuminating example of such an
issue.

Radiation is one of the most important cooling mechanisms during
switching, as it redistributes the heat produced by the current
flowing through the plasma. Hence, accurate models of radiation are
of fundamental importance to simulate the arc behavior, which is,
due to the physics of radiation at the temperatures present in the
plasma, a formidable task. Typically, the core of the arc is heated
up to 25,000$\,$K, meaning that the electromagnetic radiation
emitted by ions, atoms, and molecules of several different species
(e.g. nitrogen, oxygen, or copper) have to be taken into account.
The relevant window of the electromagnetic spectrum ranges from
$3\,10^{13}\,$Hz-$6\,10^{15}\,$Hz, corresponding to wavelengths
between $10^{-5}\,$m and $5\,10^{-8}\,$m, i.e. from infrared to
ultraviolet.\\
The main difficulty in simulating the electromagnetic radiation
emitted by an arc derives from the complexity of the emission
spectrum, where the relevant property, the absorption coefficient,
changes by many orders of magnitude at spectral lines of which
several 10,000 exist in the range under consideration. The
propagation of radiative heat in space for each frequency is
modeled, under assumptions that are reasonable for the arcing
phenomena encountered in switchgear devices, by a first-order
differential equation taking into account the absorption and the
emission of radiation along the direction of propagation. The energy
removed from the arc is with this defined by the temperature,
pressure, and composition distribution within. Due to the complexity
of the emission spectrum, a simple discretization according to
frequencies leads to hundred thousands of very thin frequency bands;
within each one of such bands the absorption coefficient can be
assumed to be constant for fixed temperature, pressure and
composition of the gaseous medium. This however, leads to the same
number of three dimensional field equations which needed to be
solved. Given the large number of finite volumes that have to be
considered in CFD simulations of realistic geometries (see e.g.
\cite{OBEGGPSW14}), the use of such a large-scale radiation model is
not feasible. Hence, there is the need to derive models for the
radiative heat transfer with much lower complexity, possibly without
compromising too much the accuracy as compared with the large-scale
model. This issue has been tackled by several contributions in the
literature \cite{JCGB14,ReGF12,FiJa91,NoIo08}. Most of the existing
approaches consist in discretizing the fraction of the EM spectrum
of interest into few bands, and assuming for each one of them some
averaged absorption properties. The mentioned approaches have the
advantage of being quite simple to implement, however  in principle
one should optimally choose ad-hoc different bands and averaged
absorption coefficients as a function of pressure, temperature and
chemical composition of the considered medium, since the radiation
parameters are affected by all these aspects. If the bands and/or
the averaged absorption coefficients are not chosen in an
appropriate way, the model accuracy is worse and more bands are
needed to improve it, resulting in a relatively large number of
bands (6-10) in order to achieve accurate results with respect to
the original, large-scale model. Hence, this procedure can be time
consuming, not trivial to carry out in a systematic way, and
ultimately suboptimal in terms of complexity/accuracy compromise.

In this paper, we present a new approach to derive small-scale,
band-averaged models of the radiative heat transfer. We first
describe the problem of radiation modeling from a novel perspective,
where the aim is to approximate the input-output behavior of a large
scale, linear-parameter-varying (LPV) dynamical system with that of
a low-order one. The large scale system has one input (black-body
intensity), one output (radiated intensity), three scheduling
parameters (temperature, pressure and composition), and a large
number of internal states (one for each considered frequency of the
EM spectrum), while the low-order LPV system has the same input,
output and scheduling parameters, but just a handful of internal
states. From this point of view, the problem can be classified as a
model-order reduction one \cite{Moor81}. Then, using classical tools
for the analysis of signals and dynamical systems, we provide
evidence that indeed models with quite low order (typically 2-3
bands) can be already good enough to capture the main behavior of
the full-order model. Finally, we tackle the order reduction problem
by using nonlinear system identification techniques (see e.g.
\cite{Nelle01}), where we collect input-output data from the
large-scale system and use it to identify the parameters of the
reduced-order model. The approach results in a nonlinear
optimization problem (nonlinear program - NLP) with a smooth
non-convex cost function and convex constraints, which are needed to
preserve the physical consistency of the reduced-order model.
We show through examples that the obtained
reduced-order models enjoy a high accuracy with respect to the
full-order one, while greatly reducing the computational times. As
compared with the existing approaches, the method proposed here has
the significant advantage of being systematic, i.e. there is no need
to tailor or tune it for each different composition of the
absorbing/emitting medium. The main user-defined parameter is the
desired number of frequency bands in the reduced-order model, which
can be then increased gradually until the desired tradeoff in terms
of model quality vs. complexity is reached.

The paper is organized as follows. Section \ref{SS:ProbForm}
provides a description of the problem we want to solve. In Section
\ref{S:SystPersp}, such a problem is analyzed from a system's
perspective and connections are made to the  order reduction of a
large-scale LPV dynamical system. The proposed computational
approach is described in Section \ref{S:NLmodred}, finally results
are presented in Section \ref{S:results} and conclusions and future
developments are discussed in Section \ref{S:conclusions}.

\subsection{Problem formulation}\label{SS:ProbForm}

\begin{figure}[hbt]
\centerline{ \includegraphics[clip,width=13cm]{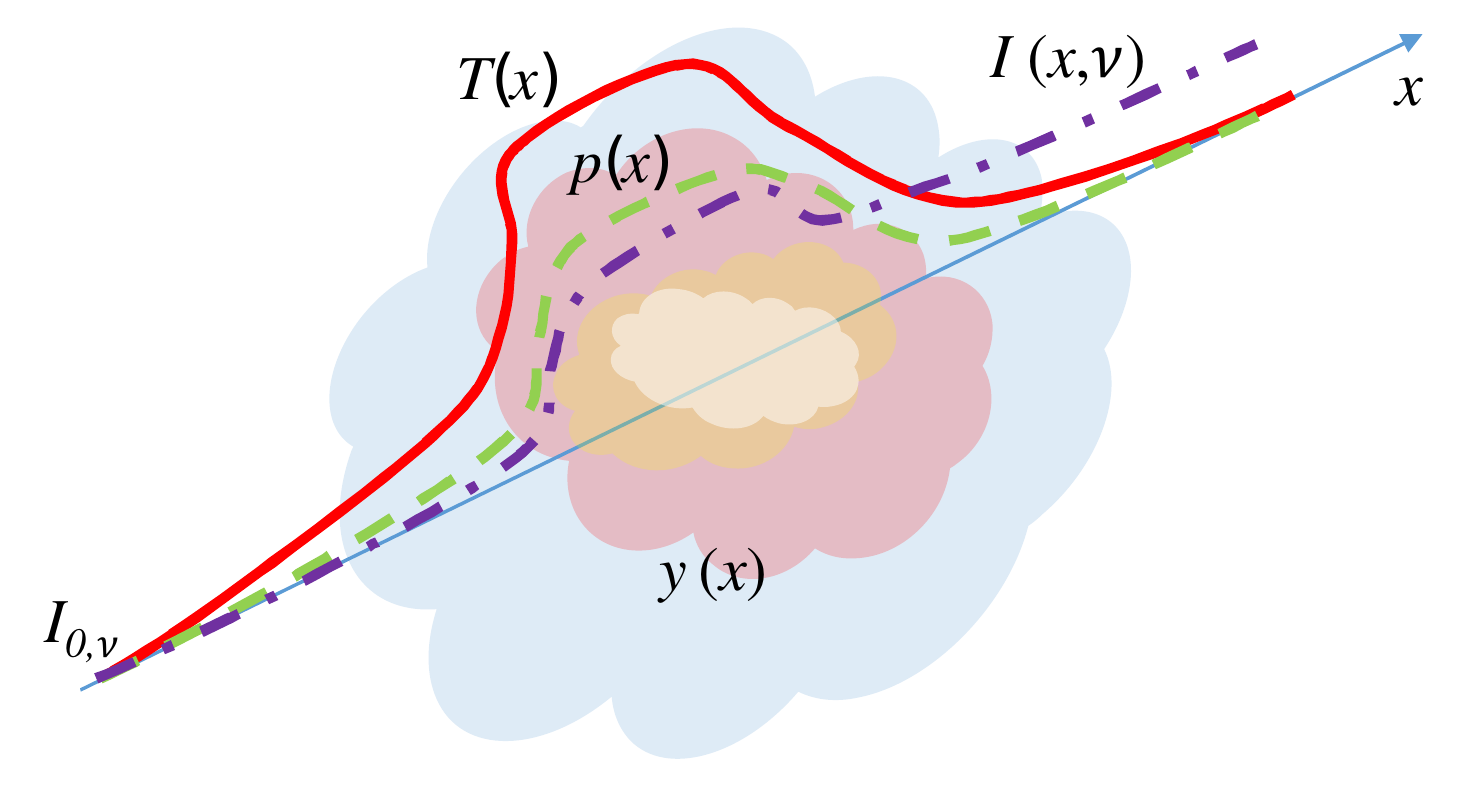} }
\protect\caption{\small Sketch of the considered problem. The aim is
to compute the spatial distribution of the radiated intensity
$I(x,\nu)$ (dash-dotted line) for a given frequency $\nu$ of the EM
spectrum through a gaseous medium along a given direction with
coordinate $x$. The temperature $T$ (solid line), pressure $p$
(dashed) and composition $y$ (balloons of different colors) of the
gas depend on $x$. The intensity at $x=0$, $I_{0,\nu}$, is a known
initial condition. \normalsize}\label{F:line_propagation}
\end{figure}

Let us consider a region in space containing a hot gaseous mixture
of $r\in\mathbb{N}$ different components. Each component is present
in a fraction (e.g. of mass or of mole)
$y_i\in[0,\,1],\,i=1,\ldots,r,$. We indicate with
$y=[y_1,\ldots,y_r]^T$ the composition vector of the medium. We
consider a line of propagation along which we want to compute the
radiated heat, and denote with $x\in\mathbb{R}$ the position of a
point lying on such a line. As the line crosses the hot matter, the
temperature $T(x)\in\mathbb{R}^+$, pressure $p(x)\in\mathbb{R}^+$
and composition $y(x)\in[0,\,1]^r$ change with $x$ (see Figure
\ref{F:line_propagation} for a graphical visualization).
The temperature, pressure and composition of the gas lie in some
sets of interest, $\mathcal{T},\,\mathcal{P}$ and $\mathcal{Y}$
respectively, defined as follows:
\begin{equation}\label{E:sets}
\begin{array}{rcl}
\mathcal{T}&=&[T_\text{min},\,T_\text{max}]\subset\mathbb{R}^+\\
\mathcal{P}&=&[p_\text{min},\,p_\text{max}]\subset\mathbb{R}^+\\
\mathcal{Y}&=&\left\{y\in[0,\,1]^r:\sum\limits_{i=1}^ry_i=1\right\}\\
\end{array}
\end{equation}
Typical values defining the sets $\mathcal{T}$ and $\mathcal{P}$ are
$T_\text{min}=300\,$K, $T_\text{max}=25,000\,$K,
$p_\text{min}=10^4\,$Pa and $p_\text{max}=10^7\,$Pa.

Assuming without loss of generality that at $x=0$ the radiated heat
intensity, for a given frequency $\nu$ of the EM spectrum, is equal
to a given value $I_{\nu,0}$, the distribution of the intensity as a
function of $x$ can be computed through the following ordinary
differential equation:
\begin{equation}\label{E:ODE_intensity}
\dfrac{dI(x,\nu)}{dx}=\alpha\left(T(x),p(x),y(x),\nu\right)\left(I_{bb}(T(x),\nu)-I(x,\nu)\right),\;\;\;I(0,\nu)=I_{\nu,0},
\end{equation}
where $\alpha(\cdot,\cdot,\cdot,\nu)$ is the absorption coefficient
and $I_{bb}(\cdot,\nu)$ is the black body intensity, both pertaining
to the frequency $\nu$. For a fixed value of $\nu$, $\alpha$ depends
on temperature, pressure and composition, while $I_{bb}$ is a
function of temperature only. 

In equation \eqref{E:ODE_intensity}, it is implicitly assumed that
the radiation distribution has already converged to a steady state.
Since the temporal dynamics of the radiated heat are much faster
than the timescale of the phenomena of interest in arc simulations,
such an assumption is reasonable. Similarly, scattering effects have
been also neglected as they represent a negligible term for the
application considered here. For a more complete form of equation
\eqref{E:ODE_intensity} the interested reader is referred to
\cite{SiHo92}.

\begin{figure}[hbt]
\centerline{ \includegraphics[clip,width=13cm]{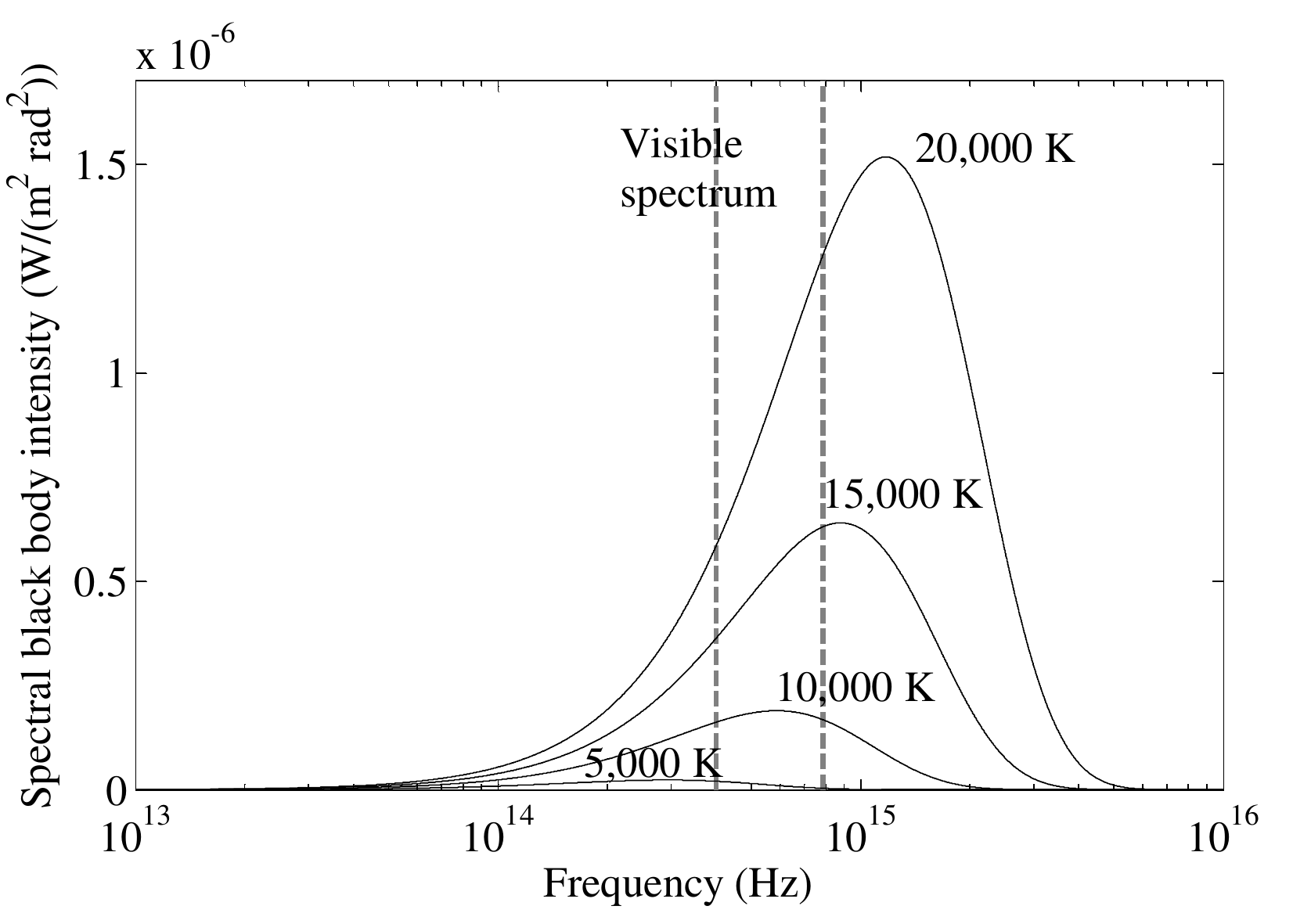} }
\protect\caption{\small Planck's law \eqref{E:Planck's law} (solid
lines) computed for different temperature values. The frequency band
of visible light is also shown (between gray dashed
lines).\normalsize}\label{Fig:Planckslaw}
\end{figure}

Thus, for each frequency $\nu$ the corresponding ODE
\eqref{E:ODE_intensity} is characterized by two parameters, namely
the black body intensity $I_{bb}$ and the absorption coefficient
$\alpha$. The dependency of the latter on $T(x),\,p(x),\,y(s)$
renders the equation nonlinear. From a physical perspective, the
black body intensity is a source of radiation, while the current
intensity $I(x,\nu)$ is a sink: the change of radiated intensity in
an infinitesimal space interval $dx$ is the difference between these
two contributions, scaled by the absorption coefficient. The
black-body intensity as a function of the frequency $\nu$ and of
temperature $T$ is given by Planck's law \cite{SiHo92}:
\begin{equation}\label{E:Planck's law}
I_{bb}(T,\nu)=\dfrac{2\,h\,\nu^3}{c^2}\dfrac{1}{e^{\frac{h\,\nu}{k_B\,T}}-1},
\end{equation}
where $h\simeq6.62\,10^{-34}\,$J$\,$s$^{-1}$ is the Planck constant,
$c\simeq10^8\,$m$\,$s$^{-1}$ is the speed of light in vacuum and
$k_B\simeq1.38\,10^{-23}\,$J$\,$K$^{-1}$ is the Boltzmann constant.
The behavior of $I_{bb}(T,\nu)$ as a function of $\nu$ for some
temperature values is shown in Figure \ref{Fig:Planckslaw}.
For the purpose of this study, the total black-body intensity
$\overline{I}_{bb}(T)$ per unit of surface and solid angle
(black-body radiance) is also needed, given by Stefan-Boltzmann's
law:
\begin{equation}\label{E:total_bb}
\overline{I}_{bb}(T)\doteq\int\limits_0^\infty
I_{bb}(T,\nu)d\nu=\dfrac{\sigma_{SB}}{\pi}T^4
\end{equation}
where $\sigma_{SB}\simeq5.67\,10^{-8}\,$W$\,$m$^{-2}\,$K$^{-4}$ is
the Stefan-Boltzmann constant. From equation \eqref{E:Planck's law}
and Figure \ref{Fig:Planckslaw} it can be noted that, for the
temperature values of interest for arc simulations in circuit
breakers (3,000K-25,000K), most of the black-body radiance (i.e. the
area below the solid curves in Figure \ref{Fig:Planckslaw}) is
contributed by the frequencies in the range of approximately
$3\,10^{13}\,$Hz-$6\,10^{15}\,$Hz.

\begin{figure}[hbt]
\centerline{
\includegraphics[clip,width=13cm]{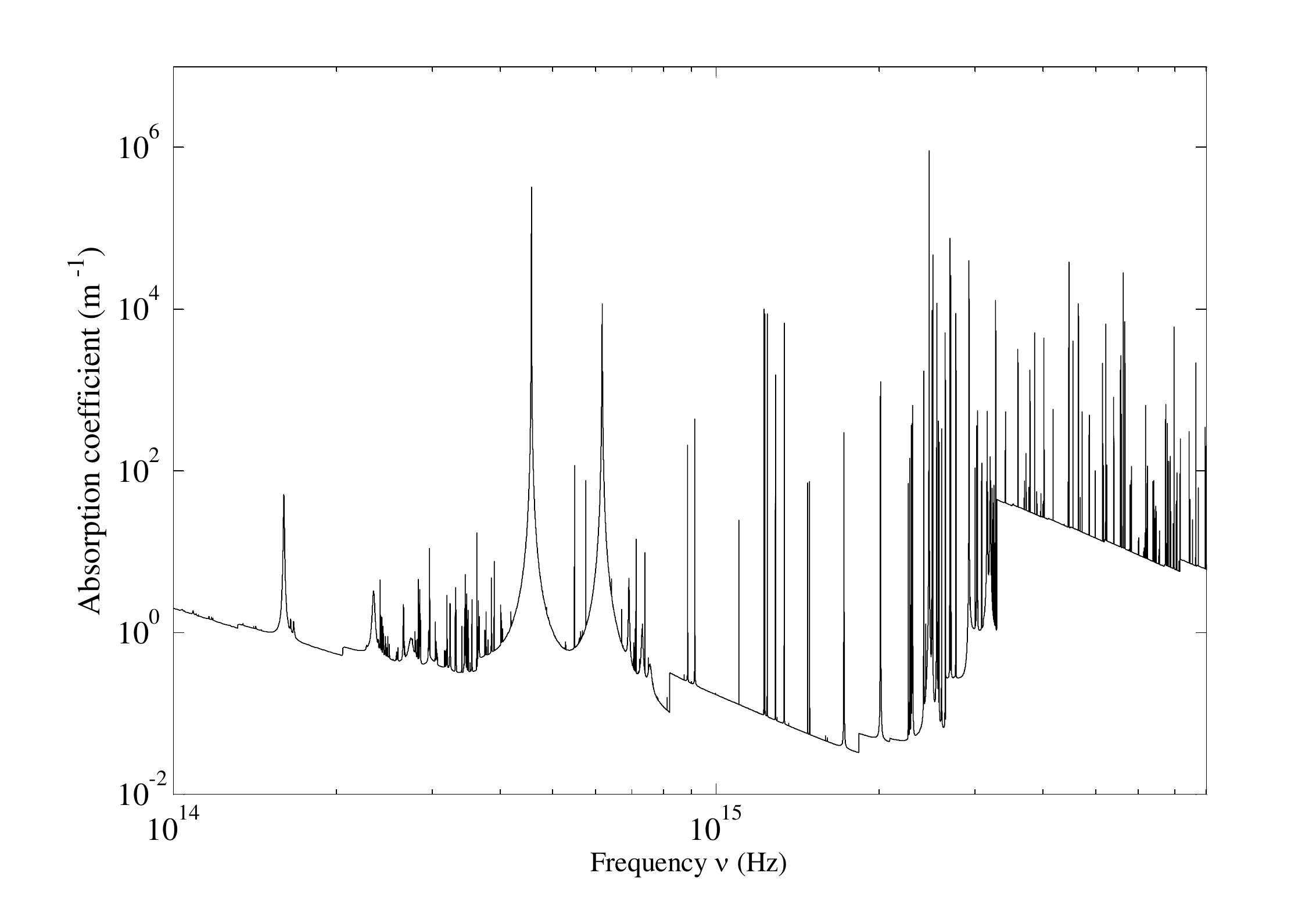} }
\protect\caption{\small Example of absorption spectrum for a mixture
of 50\% silver, 25\% air and 25\% hydrogen at 16,300 K and
10$^5\,$Pa.\normalsize}\label{Fig:abs_coeff}
\end{figure}

\begin{figure}[hbt]
\centerline{ \includegraphics[clip,width=13cm]{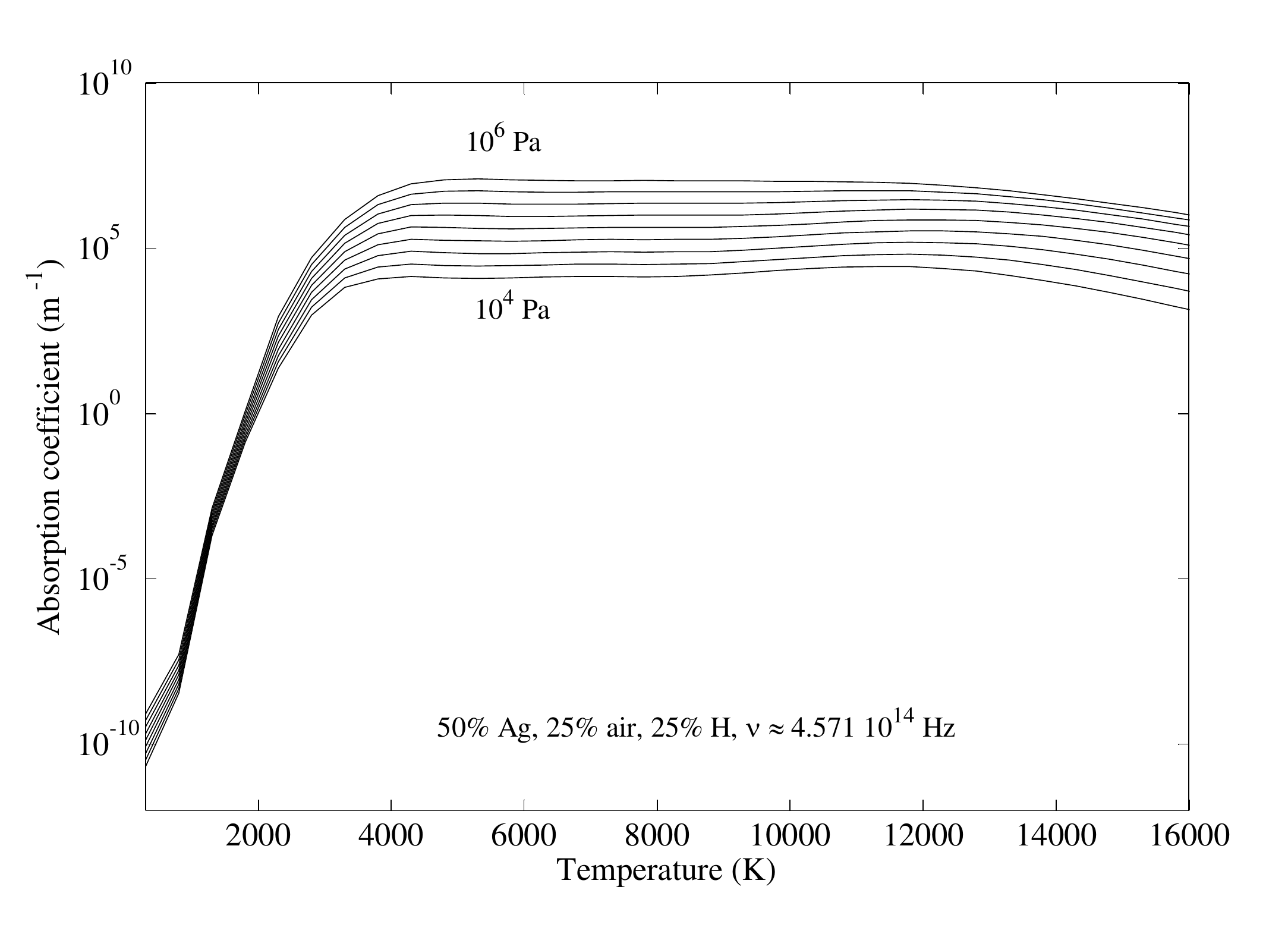} }
\protect\caption{\small Example of absorption coefficient for a
mixture of 50\% silver, 25\% air and 25\% hydrogen as a function of
temperature and pressure, for
$\nu=4.571\,10^{14}\,$Hz.\normalsize}\label{Fig:abs_coeff_T}
\end{figure}

While the dependence of the black-body intensity on $\nu$ and $T$ is
well-known from the above physical laws, the absorption coefficient
$\alpha$ is a much more uncertain quantity. Gaseous media have
absorption spectra (i.e. the function relating $\alpha$ to $\nu$)
characterized by sharp lines at specific frequencies, whose values
depend on the composition of the mixture and whose number is
typically very large. Moreover, the absorption coefficient at each
of such frequencies depends strongly on temperature and (less
markedly) on pressure. Broadening effects due to pressure are also
important as they contribute to a spread of the absorption lines
over the nearby frequencies. Physical models to capture such complex
absorption spectra have been proposed in the literature,
\cite{AuLo94,MeMa02,AuMa04,BCBVG12}. Moreover, databases of
experimentally measured absorption data for several mixtures in
specific frequency bands are available, typically for relatively low
temperatures \cite{NISTspectrum}. In the following, we will assume
that $\alpha$ is a known function of $\,T,\,p,\,y,\nu$, in the
domain of interest; we call such information the ``base data''
$\mathcal{D}$:
\begin{equation}\label{E:base_data}
\mathcal{D}\doteq\left\{\alpha(T,p,y,\nu);\;\forall\,
T\in\mathcal{T},\,
\forall\,p\in\mathcal{P},\,\forall\,y\in\mathcal{Y},\,\forall\,\nu\in\mathbb{R}^+\right\}
\end{equation}
As an example of information contained in $\mathcal{D}$, Figure
\ref{Fig:abs_coeff} shows the absorption spectrum for a mixture of
50\% silver, 25\% air and 25\% hydrogen at 16,300$\,$K and
10$^5\,$Pa. For a fixed frequency of such a spectrum, Figure
\ref{Fig:abs_coeff_T} shows the corresponding absorption coefficient
as a function of temperature and pressure (note the logarithmic
scale for the absorption coefficient in this plot).

Given the base data, the total radiated heat intensity $I_{tot}(x)$
is computed by integrating the intensity contributed by each
frequency of the EM spectrum:
\begin{equation}\label{E:total_I}
I_{tot}(x)\doteq\int\limits_0^\infty I(x,\nu)d\nu.
\end{equation}
In order to obtain a tractable problem, a first step is to adopt a
fine discretization of the frequency domain such that
\eqref{E:total_I} becomes a sum over a finite number of terms. Since
the radiated intensity is relevant only in a specific region of the
spectrum (compare Figure \ref{Fig:Planckslaw}), the frequency
discretization can be coarser outside such a region and finer
inside, in particular around the peaks of the absorption spectrum.
In this way, a finite number $N$ of frequency values
$\nu_i,\,i=1,\ldots,N$ are considered, each one being the middle
point of an interval $\Delta\nu_i$ of the spectrum. Then, equation
\eqref{E:total_I} can be re-written as:
\begin{equation}\label{E:total_I_FOM}
I(x)\simeq\sum\limits_{i=0}^N I(x,\nu_i)\Delta\nu_i.
\end{equation}
Equation \eqref{E:total_I_FOM}, together with
\eqref{E:ODE_intensity}-\eqref{E:Planck's law} and with the base
data $\mathcal{D}$ \eqref{E:base_data} evaluated at
$\nu_i,i=1,\ldots,N$, form a high-dimensional model of the radiated
heat intensity, named the Full-Order Model (FOM). The high
dimensionality of the FOM comes from the fact that a large number
$N$ of frequencies is taken into account in the discretization, so
that the approximation error is small. Typical values of $N$ for the
conditions of our interest are in the order of 1-2$\,10^5$. Due to
the large number of considered frequencies, the FOM can be used
effectively only in very simple cases, for example in
one-dimensional problems. In fact, in order to solve the radiative
distribution in two- and three-dimensional cases, as it is needed
for example in CFD simulations of real devices, one would have to
solve many equations (whose nature depends on the chosen method,
like the so-called P-N methods, the discrete ordinate method or the
photohydrodynamic approach, see \cite{SiHo92}, \cite{ChKa09}) for
each one of the considered frequencies, leading quickly to an
intractable problem. In this paper, we tackle this issue by deriving
a method to compute low complexity models for the radiative heat
transfer. More specifically, we consider the following problem:

\begin{problem}\label{P:problem}: Given the full-order model defined by
\eqref{E:ODE_intensity}-\eqref{E:Planck's law},
\eqref{E:total_I_FOM} and the base data $\mathcal{D}$, derive a
model with the same structure, i.e. where the total intensity is the
sum of a finite number of contributions obtained by partitioning
 the frequency domain, but where the number of such
partitions is very small, while still capturing accurately the total
radiated heat intensity.
\end{problem}

We call such a simplified model the Reduced Order Model (ROM). The
ROM can then be effectively used to model the propagation of heat
via radiation in many applications, including full three-dimensional
simulations of plasma arcs encountered in switchgear devices.

We remark that, in light of Problem \ref{P:problem}, in this work we
will consider the radiated intensity computed with the FOM as
``exact''. We will thus evaluate the quality of a given ROM by
assessing the discrepancy between the radiative heat intensity given
by the latter and the one given by the FOM. In other words, we don't
consider here the accuracy of the FOM with respect to the real-world
behavior of gaseous media. Indeed, the topic of modeling accurately
and/or measuring the absorption spectrum of a given gas as a
function of temperature and pressure (i.e. to compute the base data
$\mathcal{D}$) is by itself an important and active research area
\cite{AuLo94,MeMa02,AuMa04,BCBVG12}, however it is outside the scope
of this work, which is focused on the simplification of the FOM into
a computationally tractable model. On the other hand, the method
presented in this paper does not depend on the specific absorption
spectrum, i.e. it can be applied systematically to any base data,
and it yields quite small discrepancies between the derived ROM and
the employed FOM, such that the error between the ROM and the real
behavior of the considered medium depends ultimately only on the
quality of the FOM.

\section{A system's perspective of radiative heat transfer}\label{S:SystPersp}

\subsection{Equivalent input-output models of the radiative heat transfer}\label{S:LPV}
\begin{figure}[hbt]
\centerline{ \includegraphics[clip,width=13cm]{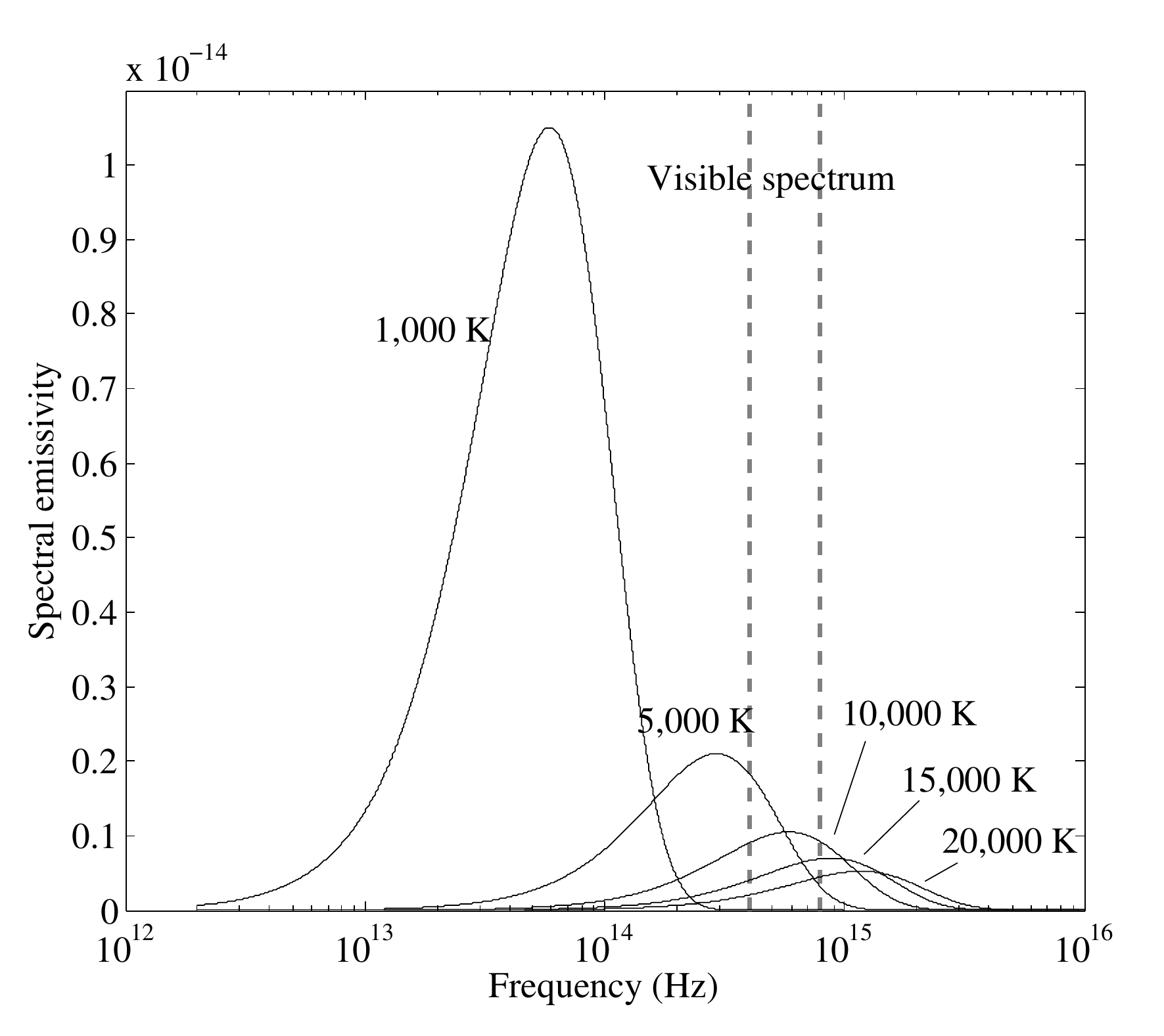} }
\protect\caption{\small Spectral emissivity $e(T,\nu)$ as a function
of EM frequency (solid lines), computed for different temperature
values. The frequency band of visible light is also shown (between
gray dashed lines).\normalsize}\label{F:spec_emiss}
\end{figure}

\begin{figure}[hbt]
\centerline{ \includegraphics[clip,width=14.5cm]{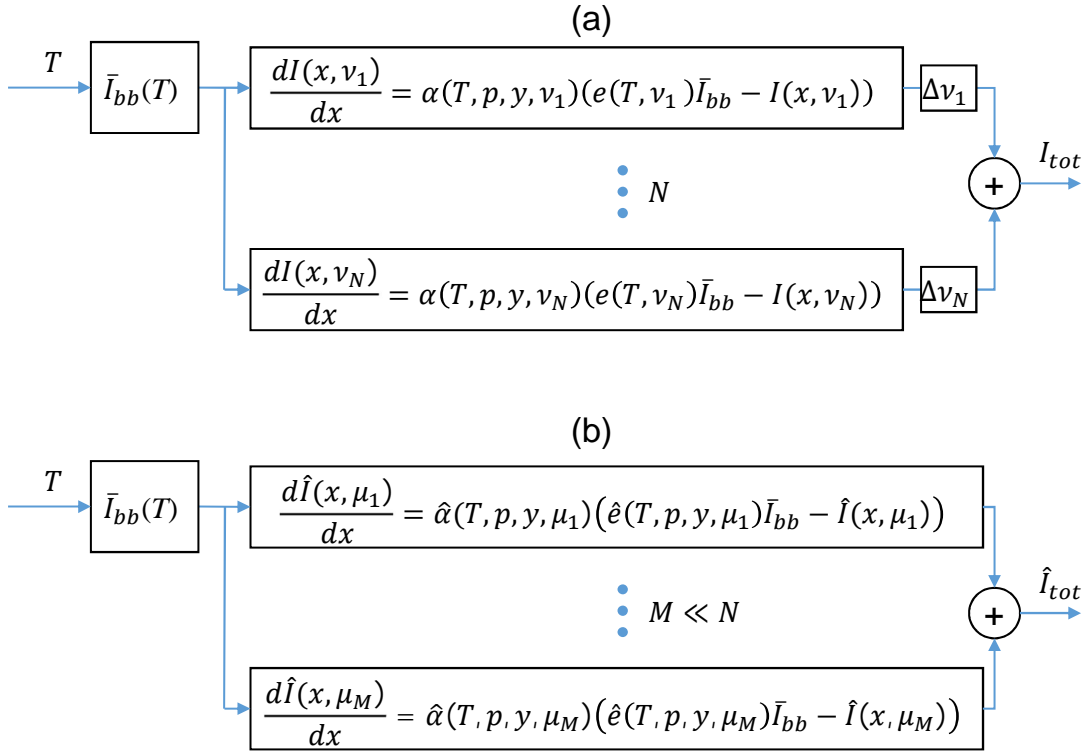}
} \protect\caption{\small Equivalent block diagram of the (a)
full-order model, FOM, and (b) reduced-order model, ROM. The
dependence of $T,p,y$ on $x$ is omitted for the sake of
readability.\normalsize}\label{F:FOM_block_diagram}
\end{figure}

As a preliminary step to address Problem \ref{P:problem}, we
re-write the FOM in a slightly different form, which is convenient
to show that this model can be seen as a single-input, single-output
(SISO) Linear Parameter-Varying dynamical system. Let us define the
spectral emissivity $e(T,\nu)$ as:
\begin{equation}\label{E:spec_emis}
e(T,\nu)\doteq\dfrac{I_{bb}(T,\nu)}{\overline{I}_{bb}(T)},
\end{equation}
i.e. the ratio between the total black body intensity and the one
pertaining to each frequency of the EM spectrum. By construction we
have $e(T,\nu)\in(0,1),\;\forall (T,\nu)$ and $\int_\nu
e(T,\nu)d\nu=1\,\forall T$. The course of $e(T,\nu)$ is illustrated
in Figure \ref{F:spec_emiss}.\\
By inserting the spectral emissivity \eqref{E:spec_emis} in equation
\eqref{E:ODE_intensity}, the radiated heat intensity for each
frequency value $\nu_i,i=1,\ldots,N$ of the FOM can be equivalently
computed as:
\begin{equation}\label{E:ODE_intensity_spec_emis}
\dfrac{dI(x,\nu_i)}{dx}=\alpha\left(T(x),p(x),y(x),\nu_i\right)\left(e(T(x),\nu_i)\overline{I}_{bb}(T(x))-I(x,\nu_i)\right),\;I(0,\nu_i)=I_{\nu_i,0}.
\end{equation}
Equation \eqref{E:ODE_intensity_spec_emis} is just a re-writing of
\eqref{E:ODE_intensity} but, together with equation
\eqref{E:total_I_FOM}, it highlights the fact that the FOM can be
seen as a dynamical system where $x$ is the independent variable,
$T(x)$ is an exogenous scalar input, $I(x,\nu_i),i=1,\ldots,N$ are
$N$ internal states, $I_{tot}(x)$ a scalar output, and $p(x),y(x)$
are space-dependent parameters. In virtue of the fact that the
function relating the temperature $T$ to the total black-body
intensity $\overline{I}_{bb}$ is known, one can consider the latter
as input to the system hence putting into evidence the Linear
Parameter Varying (LPV) structure of the model:
\begin{equation}\label{E:FOM_structure}
\begin{array}{rcll}
\dfrac{d\overline{I}(x)}{dx}&=&A(T(x),p(x),y(x))\;\overline{I}(x)+B(T(x),p(x),y(x))\;\overline{I}_{bb}(T(x))\\
I_{tot}&=&C\;\overline{I}(x),
\end{array}
\end{equation}
where
\begin{equation}\label{E:FOM_matrices}
\begin{array}{rcl}
\overline{I}(x)&\doteq&\left[I(x,\nu_1),\ldots,I(x,\nu_N)\right]^T\in\mathbb{R}^{N\times1}\\
A(T,p,y)&\doteq&\text{diag}\left(\left[-\alpha(T,p,y,\nu_1),\ldots,-\alpha(T,p,y,\nu_N)\right]^T\right)\in\mathbb{R}^{N\times
 N}\\
B(T,p,y)&\doteq&\left[\alpha(T,p,y,\nu_1)e(T,\nu_1),\ldots,\alpha(T,p,y,\nu_N)e(T,\nu_N)\right]^T\in\mathbb{R}^{N\times1}\\
C&\doteq&\left[\Delta\nu_1,\ldots,\Delta\nu_N\right]\in\mathbb{R}^{1\times
N}.
\end{array}
\end{equation}
A block-diagram of the FOM from this new point of view is shown in
Figure \ref{F:FOM_block_diagram}(a).
Such a system's perspective of the radiative heat transfer equations
is the fundamental step at the basis of all the developments
described in the remainder of this paper. We can now introduce more
precisely the structure of a candidate reduced-order model (ROM)
meant to approximate the FOM. Let us denote the \emph{bands} of the
ROM with $\mu_1,\ldots,\mu_M$, with $M\ll N$ (e.g. $M=2$), and the
corresponding vector of inner states with
$\hat{\overline{I}}(x)\doteq[\hat{I}(x,\mu_1),\ldots,\hat{I}(x,\mu_M)]^T\in\mathbb{R}^{M\times1}$.
 Then, we can
write the equations describing the ROM as:
\begin{equation}\label{E:ROM_structure}
\begin{array}{rcll}
\dfrac{d\hat{\overline{I}}(x)}{dx}&=&\hat{A}(T(x),p(x),y(x))\;\hat{\overline{I}}(x)+\hat{B}(T(x),p(x),y(x))\;\overline{I}_{bb}(T(x))\\
\hat{I}_{tot}(x)&=&\hat{C}\;\hat{\overline{I}}(x),
\end{array}
\end{equation}
where
\begin{equation}\label{E:ROM_matrices}
\begin{array}{rcl}
\hat{A}(T,p,y)&\doteq&\text{diag}\left(\left[-\hat{\alpha}(T,p,y,\mu_1),\ldots,-\hat{\alpha}(T,p,y,\mu_N)\right]^T\right)\in\mathbb{R}^{M\times
 M}\\
\hat{B}(T,p,y)&\doteq&\left[\hat{\alpha}(T,p,y,\mu_1)\hat{e}(T,p,y,\mu_1),\ldots,\hat{\alpha}(T,p,y,\mu_N)\hat{e}(T,p,y,\nu_N)\right]^T\in\mathbb{R}^{M\times1}\\
\hat{C}&\doteq&\left[1,\ldots,1\right]\in\mathbb{R}^{1\times M}.
\end{array}
\end{equation}
Figure \ref{F:FOM_block_diagram}(b) gives a graphical representation
of the ROM. Practically speaking, from
\eqref{E:ROM_structure}-\eqref{E:ROM_matrices} one can see that each
component of vector $\hat{\overline{I}}(x)$ accounts for a certain
portion of the total radiated heat intensity, in complete analogy
with the FOM \eqref{E:FOM_structure}, the only difference being the
number of frequency bands, which in the ROM is much smaller than in
the FOM.
 For each band $\mu_i$, the parameters
$\hat{\alpha}(T,p,y,\mu_i)$ and $\hat{e}(T,p,y,\mu_i)$ have thus the
meaning of ``equivalent'' absorption coefficient and spectral
emissivity, respectively. The task of deriving a ROM for the
radiated heat transfer from the FOM is referred to as \emph{model
order reduction}. Computing a ROM is equivalent to assigning
suitable values to $\hat{\alpha}$ and $\hat{e}$ as a function of the
underlying parameters $T,p$ and $y$. In particular the collection of
all bands has to form a non-overlapping partition covering the whole
EM spectrum. Then, the value of $\hat{e}(T,p,y,\mu_i)$ has to
correspond to the integral of the spectral black body intensity
\eqref{E:Planck's law} over the frequency band pertaining to
$\mu_i$. Hence, one can equivalently state that computing a ROM
amounts to choose, for each pair of pressure and composition values,
a partition of the EM spectrum (which defines the equivalent
emissivity as a function of temperature) and the courses of the
corresponding equivalent absorption coefficients as a function of
temperature. In previous contributions in the literature, e.g.
\cite{NoIo08,ReGF12}, the task of defining the ROM has been carried
out by picking a finite number of bands covering the EM spectrum and
then computing the equivalent absorption coefficients
$\hat{\alpha}(T,p,y,\mu_i)$ through some averaging procedure on the
portion of the absorption spectrum contained in each band. This
approach is simple to implement but it has the drawback of not being
systematic, since both the choice of the band cuts and the averaging
of the absorption coefficients have to be made by the user, without
an immediate link to the accuracy of the resulting ROM. In the next
sections, we will present a new order reduction approach to compute
the partitioning of the EM spectrum and the values of
$\hat{\alpha}(T,p,y,\mu)$ in a systematic way, that yields quite
accurate results as compared with the FOM.

About this last point, i.e. the accuracy of the ROM, we note that
the input of the FOM and of the ROM is exactly the same,
corresponding to the total black-body radiation for the considered
temperature profile, $\overline{I}_{bb}(T(x))$. Therefore, it is
quite intuitive that the discrepancy between the total intensity
given by the FOM, $I_{tot}(x)$, and the one predicted by the ROM,
$\hat{I}_{tot}(x)$, for the same temperature profile $T(x)$ (i.e.
the same distribution of $\overline{I}_{bb}(T(x))$), represents a
reasonable indicator of the accuracy of the reduced order model. In
other words, the error signal $\Delta I(x)\doteq
I_{tot}(x)-\hat{I}_{tot}(x)$ will be considered to evaluate the
goodness of a given ROM. This choice is motivated by the fact that
the total radiated heat intensity is the main quantity of interest
predicted by the ROM when it is embedded in multi-physics
simulations of arc plasma, since it is used to compute the heat
transferred from the plasma volume to the walls, and also (through
its divergence) the heat redistributed within the plasma volume.
Hence, the ROM should reproduce this quantity as accurately as
possible with respect to the FOM, given the same spatial
distribution of temperature, pressure and chemical composition.

Before going to the details of the proposed method to derive the
ROM, a sensible question to be addressed is whether the
approximation problem we are dealing with has a reasonably good
solution or not. More specifically, recall that we aim to
approximate the input-output behavior of a large scale system, with
hundred of thousands of internal states, with that of a small-scale
one, with at most a handful of internal states. It is not
immediately clear if there exist such a low-order ROM still capable
of delivering high approximation accuracy, since this aspect depends
on the characteristics of the FOM, i.e. on the underlying physics of
the radiative heat transfer. In the next section, we exploit the
system's perspective described above to provide an intuition that
indeed a ROM with a handful of bands can capture most of the
input-output behavior of the FOM.

\subsection{Frequency-domain analysis}\label{S:FreqDom}

Analyzing the complexity of the FOM in its form
\eqref{E:FOM_structure} in the domain of position $x$ is not
straightforward due to the very large number of internal states,
each one following its own dynamic evolution. Besides noticing that
the FOM is given by the sum of a large number of non-interacting,
asymptotically stable first-order systems, all driven by the same
input and whose (position dependent) poles are given by
$-\alpha(T,p,y)$, little less can be said.

\begin{figure}[hbt]
\centerline{ \includegraphics[clip,width=13cm]{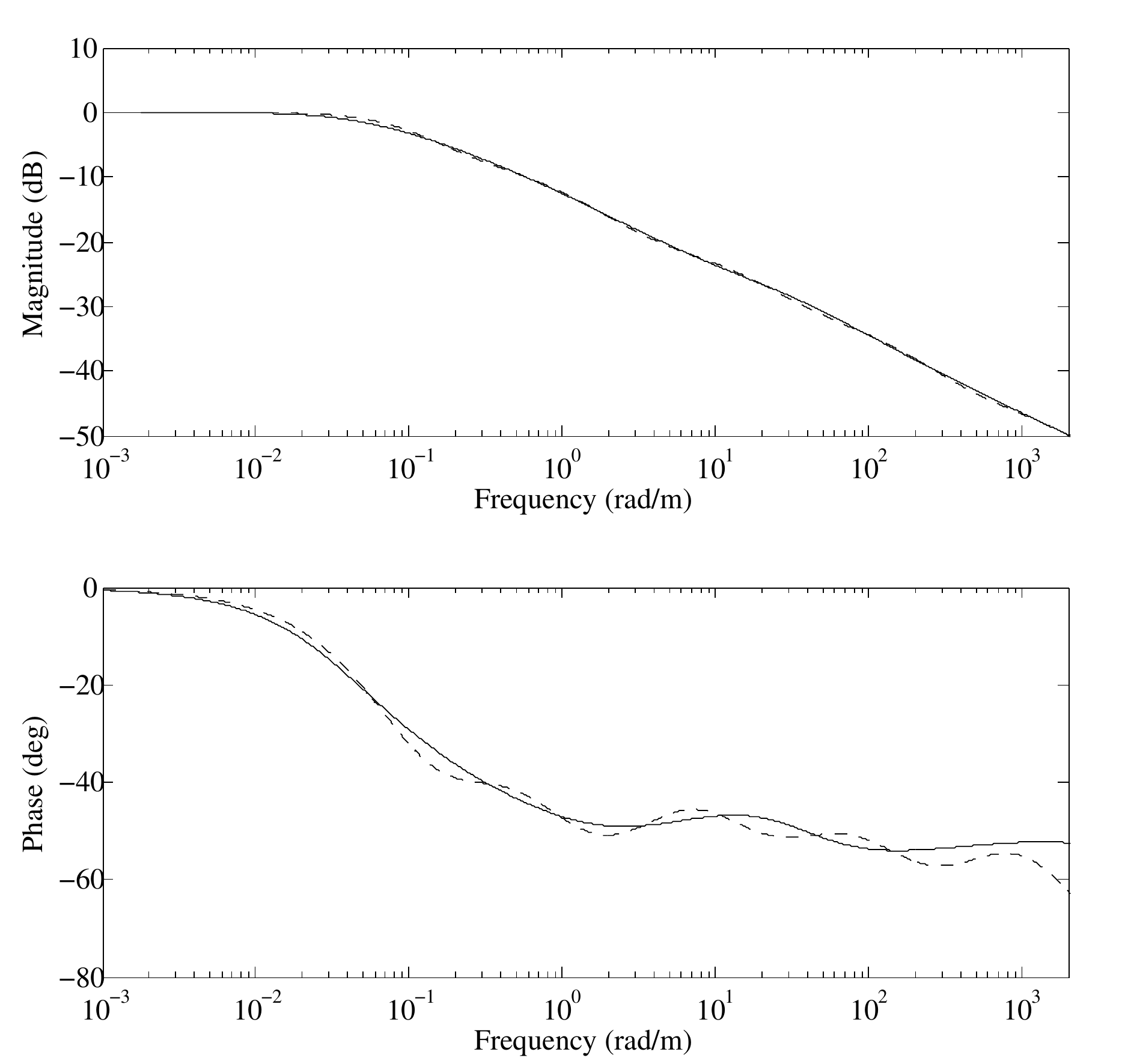}
} \protect\caption{\small Example of frequency response for a
mixture of 50\% silver, 25\% air and 25\% hydrogen at 16,300 K and
$10^5\,$Pa. Solid: full-order model; dashed: reduced-order model
with $M=5$ bands.\normalsize}\label{F:bode_FOM}
\end{figure}

However, if we consider fixed values of temperature, pressure and
composition, $\overline{T},\overline{p}$ and $\overline{y}$
respectively, and we assume that only infinitesimal perturbations of
temperature take place in space, such that the absorption properties
of the medium can be assumed constant, we immediately notice that
the FOM becomes a linear-parameter-invariant (LPI) system, to which
well-assessed tools in systems theory and signal processing can be
applied. In particular, after establishing the analogy between the
position $x$ in the FOM with the continuous time variable in
dynamical systems, we can study the input-output response of the FOM
to such infinitesimal temperature variations by applying the Laplace
transform \cite{Schi99} to its equations and deriving the transfer
function $G(s)$ from its input $I_{bb}(s)$ to its output
$I_{tot}(s)$, where $s$ is the Laplace variable:
\begin{equation}\label{E:FOM_TF}
G(s)\doteq\dfrac{I_{tot}(s)}{\overline{I}_{bb}(s)}=C\left(s\textbf{I}-A(\overline{T},\overline{p},\overline{y})\right)^{-1}B(\overline{T},\overline{p},\overline{y}),
\end{equation}
In \eqref{E:FOM_TF}, $A,B$ are the matrices given in equation
\eqref{E:FOM_matrices} and evaluated at the chosen temperature,
pressure and composition values, $\textbf{I}$ is the identity matrix
of suitable order and $\left(s\textbf{I}-A\right)^{-1}$ denotes a
matrix inverse operation. A tool commonly used  to analyze the
transfer function of a dynamical system is the Bode diagram of the
corresponding frequency response, obtained by evaluating the
magnitude and phase of $G(j\omega)$, where
$\omega=\frac{2\pi}{\tau}$ assumes here the physical meaning of the
frequency of purely sinusoidal oscillations (in space) of the
black-body intensity (i.e. of temperature), with infinitesimal
amplitude,  with the period $\tau$ measured in $m$ (e.g.
$\omega=628\,$rad/m corresponds to a period of oscillation of the
input of roughly $10^{-2}\,$m). As an example, the Bode diagram of
the FOM frequency response obtained by fixing
$\overline{T}=16,300\,$K, $\overline{p}=5\,10^5\,$Pa and
$\overline{y}$ containing 50\% silver, 25\% air and 25\% hydrogen
(whose absorption spectrum is shown in Figure \ref{Fig:abs_coeff})
is shown in Figure \ref{F:bode_FOM}.
It can be noted that the overall behavior of the FOM is that of a
low-pass filter, whose frequency response is shaped by the
contributions of the large number of internal states present in the
system. As a quantitative example regarding this case, a sinusoidal
oscillation of the black-body intensity (due to an oscillation in
temperature) with a period of 2$\pi\,10^3\,$m would give rise, in
the gaseous medium, to a sinusoidal distribution of radiated heat
intensity with unchanged amplitude and phase with respect to what
would happen in vacuum (compare Figure \ref{F:bode_FOM} with
$\omega=10^{-3}\,$rad/m), i.e. the amplitude of the oscillations
would be equal to that of the input black-body intensity, and the
spatial distribution would show almost zero phase shift. On the
other hand, the spatial distribution of radiated heat induced by a
sinusoidal oscillation of the black-body intensity with a period of
2$\pi\;$m (compare again Figure \ref{F:bode_FOM}, with
$\omega=1\,$rad/m) would have an amplitude equal to only about
$23\%$ with respect to that of the input, with a phase lag of about
$45\,$deg.

If we consider now the order reduction of such a dynamical system,
we see that this is a standard problem of model order reduction of
LPI systems, for which a well-established literature exist, see e.g.
\cite{Moor81,SaFr87}. Thus, we can use one of the existing
approaches to derive a reduced order model that approximates the FOM
for the chosen values of $\overline{T},\overline{p}$ and
$\overline{y}$. After deriving the ROM, the related transfer
function can be computed as (compare equation
\eqref{E:ROM_matrices}):
\begin{equation}\label{E:ROM_TF}
\hat{G}(s)\doteq\dfrac{\hat{I}_{tot}(s)}{\overline{I}_{bb}(s)}=\hat{C}\left(sI-\hat{A}(\overline{T},\overline{p},\overline{y})\right)^{-1}\hat{B}(\overline{T},\overline{p},\overline{y}).
\end{equation}
A comparison between the Bode diagram of the latter and that of the
FOM reveals that up to very small gain values of about $10^{-4}$,
such that the corresponding radiated heat intensity is negligible,
the frequency responses of the two models are practically
super-imposed, hence showing a very good agreement between them.
Moreover, applying this procedure for many values of
$\overline{T},\overline{p}$ and $\overline{y}$, chosen by gridding
their respective domains $\mathcal{T},\,\mathcal{P}$ and
$\mathcal{Y}$, shows that such a good agreement is obtained always
with no more than four-five bands in the ROM. A good agreement up to
a gain of about $5\%$ is obtained with just two bands in the ROM,
for oscillation periods of fractions of millimeters. Indeed, this
level of accuracy would be enough for the sake of arc plasma CFD
simulations, where the resolution of the spatial discretization of
the considered volume is of the order of $10^{-3}\,$m. An example of
the obtained results is depicted in Figure \ref{F:bode_FOM}, too. In
particular, the ROM whose frequency response is shown in the figure
has five bands, $\mu_1,\ldots,\mu_5$, and the corresponding
absorption coefficients and emissivities (and frequency boundaries)
are equal to
\[
\begin{array}{rclcrcl}
\hat{\alpha}(\overline{T},\overline{p},\overline{y},\mu_1)&=&1.1\,10^{-1}\;\text{m}^{-1}&\;&\hat{e}(\overline{T},\overline{p},\overline{y},\mu_1)&=&7.2\;10^{-1}\;\;[0,\;1.6\,10^{15}]\,\text{Hz}\\
\hat{\alpha}(\overline{T},\overline{p},\overline{y},\mu_2)&=&9.2\,10^{-1}\;\text{m}^{-1}&\;&\hat{e}(\overline{T},\overline{p},\overline{y},\mu_2)&=&2.1\;10^{-1}\;\;(1.6\,10^{15},\;2.4\,10^{15}]\,\text{Hz}\\
\hat{\alpha}(\overline{T},\overline{p},\overline{y},\mu_3)&=&1.3\,10^{1}\;\text{m}^{-1}&\;&\hat{e}(\overline{T},\overline{p},\overline{y},\mu_3)&=&5.4\;10^{-2}\;\;(2.4\,10^{15},\;3.1\,10^{15}]\,\text{Hz}\\
\hat{\alpha}(\overline{T},\overline{p},\overline{y},\mu_4)&=&1.2\,10^{2}\;\text{m}^{-1}&\;&\hat{e}(\overline{T},\overline{p},\overline{y},\mu_4)&=&1.2\;10^{-2}\;\;(3.1\,10^{15},\;3.8\,10^{15}]\,\text{Hz}\\
\hat{\alpha}(\overline{T},\overline{p},\overline{y},\mu_5)&=&1.7\,10^{3}\;\text{m}^{-1}&\;&\hat{e}(\overline{T},\overline{p},\overline{y},\mu_5)&=&3.9\;10^{-3}\;\;(3.8\,10^{15},\;+\infty)\,\text{Hz}\\
\end{array}
\]
A comparison between these values and Figure \ref{F:bode_FOM} shows
that the values of the ROM absorption coefficients correspond to the
dominant poles of the FOM, as the intuition would suggest.

Overall, the analysis reported so far provides if not a rigorous
proof at least an indication that the problem we are dealing with
has a reasonable solution, using ROMs of quite low order. In section
\ref{S:NLmodred}, we describe in details the solution approach that
we propose to deal with the linear-parameter-varying case.

\subsection{Position discretization}\label{SS:space_discr}
Before proceeding further, it is convenient to introduce the
discretized versions of the FOM and of the ROM, where the position
variable $x$ is taken at nodes $x_l,\,l\in\mathbb{N}$, that are
equally spaced by an interval $\Delta x$. The latter has to be
chosen according to the features of the problem at hand, trading off
computational speed with a sufficiently fine discretization, which
can capture well the fastest transients of the model's input and
output. As a rule of thumb, $\Delta x\simeq \delta/10$ can be
chosen, where $\delta$ is the smallest space-scale of interest in
the problem (typically in our case $\delta\simeq10^{-3}\,$m, as
mentioned above). Hence, all the space-dependent variables (i.e.
$T,p,y$) are now evaluated at discrete position values. The
discretization of the radiation models is carried out by assuming
that such variables are constant between two subsequent position
nodes,  $x_l$ and $x_{l+1}=x_l+\Delta x$, and then computing the
explicit integration of \eqref{E:FOM_structure} (for the FOM) and
\eqref{E:ROM_structure} (for the ROM). In particular, for the
full-order model we have:
\begin{equation}\label{E:FOM_structure_dt}
\begin{array}{rcll}
\overline{I}(x_{l+1})&=&A_d(T(x_l),p(x_l),y(x_l))\;\overline{I}(x_l)+B_d(T(x_l),p(x_l),y(x_l))\;\overline{I}_{bb}(T(x_l))\\
I_{tot}(x_{l})&=&C\;\overline{I}(x_{l})
\end{array}
\end{equation}
where the matrices $A_d$ and $B_d$ are computed as:
\begin{equation}\label{E:FOM_matrices_dt}
\begin{array}{rcl}
A_d(T,p,y)&\doteq&\text{diag}\left(\left[a(T,p,y,\nu_1),\ldots,a(T,p,y,\nu_N)\right]^T\right)\in\mathbb{R}^{N\times
 N}\\
\hat{B}_d(T,p,y)&\doteq&\left[(1-a(T,p,y,\nu_1))e(T,\nu_1),\ldots,(1-a(T,p,y,\nu_N))e(T,\nu_N)\right]^T\in\mathbb{R}^{N\times1}
\end{array},
\end{equation}
and
\[
a(T,p,y,\nu_i)= e^{-\alpha(T,p,y,\nu_i)\Delta x},\,i=1,\ldots,M.
\]
The matrix $C$ in \eqref{E:FOM_structure_dt} is the same as in
\eqref{E:FOM_structure}, since the output equation is static and
thus it is not changed by the discretization of the position $x$.
Similarly, for the reduced-order model we have:
\begin{equation}\label{E:ROM_structure_dt}
\begin{array}{rcll}
\hat{\overline{I}}(x_{l+1})&=&\hat{A}_d(T(x_l),p(x_l),y(x_l))\;\hat{\overline{I}}(x_l)+\hat{B}_d(T(x_l),p(x_l),y(x_l))\;\overline{I}_{bb}(T(x_l))\\
\hat{I}_{tot}(x_{l})&=&\hat{C}\;\hat{\overline{I}}(x_{l})
\end{array}
\end{equation}
where the matrices $\hat{A}_d$ and $\hat{B}_d$ are computed as:
\begin{equation}\label{E:ROM_matrices_dt}
\begin{array}{rcl}
\hat{A}_d(T,p,y)&\doteq&\text{diag}\left(\left[\hat{a}(T,p,y,\mu_1),\ldots,\hat{a}(T,p,y,\mu_M)\right]^T\right)\in\mathbb{R}^{M\times
 M}\\
\hat{B}_d(T,p,y)&\doteq&\left[(1-\hat{a}(T,p,y,\mu_1))\hat{e}(T,p,y,\mu_1),\ldots,(1-\hat{a}(T,p,y,\mu_M))\hat{e}(T,p,y,\mu_M)\right]^T\in\mathbb{R}^{M\times1}
\end{array},
\end{equation}
and
\begin{equation}\label{E:coeff_a_dt}
\hat{a}(T,p,y,\mu_i)= e^{-\hat{\alpha}(T,p,y,\mu_i)\Delta
x},\,i=1,\ldots,M.
\end{equation}
Also for the ROM the matrix $\hat{C}$ in \eqref{E:ROM_structure_dt}
is the same as in \eqref{E:ROM_structure}.

The reason why we employ such a space discretization is twofold: on
the one hand, it is needed to obtain a finite-dimensional
computational problem, on the other hand it improves the
computational efficiency of the method (in particular by using a
constant discretization step $\Delta x$).

\section{Model order reduction of the radiative heat transfer equation}\label{S:NLmodred}

\subsection{Solution approach}\label{SS:Solution}
Considering the analysis of Section \ref{S:FreqDom}, one can be
tempted to derive the ROM, i.e. the functions $\hat{\alpha}(T,p,y)$
and $\hat{e}(T,p,y)$, by gridding the domains
$\mathcal{T},\mathcal{P}$ and $\mathcal{D}$ and for each triplet
$(\overline{T},\overline{p},\overline{y})$ compute the equivalent
absorption coefficient and emissivity of the corresponding LPI
model, using well-assessed and efficient model order reduction
techniques. Then, the ROM can be obtained by interpolating among the
computed reduced-order LPI models. This approach could work well if
only pressure and composition dependence were considered, since the
sensitivity of the base data on these values is mild and hence one
can be confident that interpolating among the bands computed at
different triplets yields correct results. In other words, for a
given temperature $\overline{T}$, the absorption spectrum of the FOM
does not change dramatically between two neighboring pairs
$(\overline{p}_1,\overline{y}_1)$ and
$(\overline{p}_2,\overline{y}_2)$ within the pressure and
composition intervals of interest for the application considered
here, so that for each band $\mu_i$ of the ROM it is safe to
interpolate between the coefficients
$\hat{\alpha}(\overline{T},\overline{p}_1,\overline{y}_1,\mu_i)$ and
$\hat{\alpha}(\overline{T},\overline{p}_2,\overline{y}_2,\mu_i)$,
computed independently by means of LPI model order reduction, to
obtain the ROM coefficients for generic values of
$(\overline{T},p,y)$ with
$p\in[\overline{p}_1,\,\overline{p}_2],\;y\in[\overline{y}_1,\,\overline{y}_2]$.
However, this procedure would not achieve good results when
temperature dependence is taken into account as well. In fact, the
dependence of the absorption spectrum on temperature is very strong
(compare Figure \ref{Fig:abs_coeff_T}), so that the ROM coefficients
pertaining to the same band (e.g. $\mu_1$) but computed at two
different temperature values, even with the same pressure and
composition values, might be completely unrelated to each other, and
interpolating between them can give highly inaccurate results.
Driven by these considerations, we adopt a hybrid strategy, where we
grid the domains $\mathcal{P}$ and $\mathcal{Y}$ and for each pair
$(\overline{p},\overline{y})$ we derive the partition of the EM
spectrum in $M$ bands and the related functions
$\hat{\alpha}(T,\overline{p},\overline{y},\mu_i),\,i=1,\ldots,M$
that define the ROM. Since the corresponding FOM is now parameter
varying (because we let the temperature distribution change while
fixing only pressure and composition), this problem falls in the
class of model-order reduction of LPV systems, for which,
differently from the LPI case, few results exist in the literature
\cite{WoGG96,SaRa04}, and their practical applicability to systems
with $\approx10^5$ states, like the FOM in our problem, is not
straightforward. For the above reasoning, in the remainder of this
section it is assumed, unless otherwise stated, that a fixed pair
$(\overline{p},\overline{y})$ of pressure and composition values has
been chosen, and that the only space-varying variable is the
temperature $T$. The complete ROM  can be then obtained by repeating
the LPV order-reduction for all the pairs
$(\overline{p},\overline{y})$ chosen by gridding the respective
domains, and then interpolating among the obtained values to compute
the equivalent absorption coefficients and frequency bands for a
generic triplet $(T,p,y)$.
\begin{remark}\label{R:remark_const_py}
The simplification introduced by fixing pressure and composition is
made possible thanks to the above-discussed particular properties of
the absorption spectra of the considered gaseous media. Indeed, such
a simplification by itself is not required for the order-reduction
technique described in the following, which may straightforwardly be
applied also with varying pressure and composition, but at the price
of higher computational requirements. In our experience, using
constant pressure and composition in the order reduction computation
yields accurate enough results for the applications of interest.
\end{remark}

In the following, we propose to address the LPV order reduction
problem with a nonlinear system identification approach (see e.g.
\cite{Nelle01}), in which we search, within a given set of possible
ROMs, the one which is closest to the FOM according to a pre-defined
optimality criterion. In the next sections, this task is brought in
the form of a tractable optimization program. It has to be noted at
this point that, in the literature on system identification, there
exist several contributions devoted to the problem of identification
of LPV systems, see e.g. \cite{BaGi02} and the references therein.
However, such results are not applicable in our case, due to the
additional constraints that are present on the ROM, namely the need
to preserve a specific structure where the total radiated intensity
is the sum of the contributions given by the frequency bands.
Differently from such previous approaches, the one proposed here is
able to take into account these constraints, since it is
specifically tailored for the considered application.

\subsubsection{Cost function and model set}

In order to have an optimization problem that can be solved with
common numerical techniques, two main ingredients need to be
defined: the set of reduced-order models of the form
\eqref{E:ROM_structure_dt} where we carry out our search, denoted
with $\mathcal{H}$ (model set), and a cost function $J$ giving a
measure of how much a given ROM $H\in\mathcal{H}$ is close to the
FOM. The model set $\mathcal{H}$ should represent the limitations
that we want to impose on the ROM in order to account for the
physics of the problem. A ROM is fully characterized by its
parameters
$\hat{\alpha}(T,\overline{p},\overline{y},\mu_i),\,\hat{e}(T,\overline{p},\overline{y},\mu_i),\,i=1,\ldots,M$,
which in the considered settings are, for each band $\mu_i$,
functions of temperature only. To be consistent with the underlying
physical phenomena, the equivalent absorption coefficients should
 be positive (to retain an asymptotically stable model):
\begin{equation}\label{E:constraints_ct_func_alpha}
\begin{array}{l}
\hat{\alpha}(T,\overline{p},\overline{y},\mu_i)>0,\,\forall
T\in\mathcal{T},\,i=1,\ldots,M
\end{array}
\end{equation}
and the equivalent emissivities should lie in the interval $[0,\,1]$
and sum to one over all the considered temperature range, in analogy
with \eqref{E:spec_emis}, so that the black-body limit is not
violated:
\begin{equation}\label{E:constraints_ct_func_e}
\begin{array}{l}
0\leq\hat{e}(T,\overline{p},\overline{y},\mu_i)\leq 1,\,\forall
T\in\mathcal{T},\,i=1,\ldots,M\\
\sum\limits_{i=1}^M\hat{e}(T,\overline{p},\overline{y},\mu_i)=1,\,\forall
T\in\mathcal{T}
\end{array}
\end{equation}
Since equation \eqref{E:coeff_a_dt} is invertible, we can
equivalently consider the functions
$\hat{a}(T,\overline{p},\overline{y},\mu)$ to define the model set.
We select this alternative for the sake of computational efficiency,
as we will discuss more in details in section \ref{SSS:comp_eff}.
Then, the constraints \eqref{E:constraints_ct_func_alpha} can be
re-written as:
\begin{equation}\label{E:constraints_dt_func}
\begin{array}{l}
0\leq \hat{a}(T,\overline{p},\overline{y},\mu_i)<1,\,\forall
T\in\mathcal{T},\,i=1,\ldots,M
\end{array}
\end{equation}
As a final step to define $\mathcal{H}$, we choose a finite
parametrization of functions
$\hat{a}(T,\overline{p},\overline{y},\mu_i)$ and
$\hat{e}(T,\overline{p},\overline{y},\mu_i)$, in order to obtain a
finite dimensional model set (hence also a finite dimensional
optimization problem):
\begin{equation}\label{E:parametrization}
\begin{array}{l}
\hat{a}(T,\overline{p},\overline{y},\mu)=f_{\hat{a}}(T,\theta_{\hat{a}}(\overline{p},\overline{y},\mu)),\,\theta_{\hat{a}}(\overline{p},\overline{y},\mu)\in\mathbb{R}^{n_{\theta_{\hat{a}}}}\\
\hat{e}(T,\overline{p},\overline{y},\mu)=f_{\hat{e}}(T,\theta_{\hat{e}}(\overline{p},\overline{y},\mu)),\,\theta_{\hat{e}}(\overline{p},\overline{y},\mu)\in\mathbb{R}^{n_{\theta_{\hat{e}}}}
\end{array},
\end{equation}
where $f_{\hat{a}},\,f_{\hat{e}}$ are chosen, once again, to
tradeoff the flexibility of the parametrization with computational
efficiency. In particular, a convenient choice for $f_{\hat{a}}$ is
the class of piecewise affine functions of temperature, while
$f_{\hat{e}}$ is taken such that the corresponding parameters
$\theta_{\hat{e}}$ define the partitioning of the EM spectrum in a
finite number $M$ of bands. These choices have the advantage of
 yielding a convex model set (see section
\ref{SSS:parametrization} for details), hence improving the
computational efficiency and stability of the procedure. For a given
pair $(\overline{p},\overline{y})$, let us collect the parameters
$\theta_{\hat{a}}$ and $\theta_{\hat{e}}$ for all the ROM bands
$\mu_i$ in a single vector $\theta$:
\begin{equation}\label{E:param_vec}
\begin{array}{rcl}
\theta&\doteq&\left[
\begin{array}{c}
\theta_{\hat{a}}(\overline{p},\overline{y},\mu_1)\\
\vdots\\
\theta_{\hat{a}}(\overline{p},\overline{y},\mu_M)\\
\theta_{\hat{e}}(\overline{p},\overline{y},\mu_1)\\
\vdots\\
\theta_{\hat{e}}(\overline{p},\overline{y},\mu_M)
\end{array}\right]\in\mathbb{R}^{n_\theta},
\end{array}
\end{equation}
where $n_{\theta}=M(n_{\theta_{\hat{a}}}+n_{\theta_{\hat{e}}})$ is
the total number of optimization variables. Then, considering
\eqref{E:constraints_ct_func_e}-\eqref{E:param_vec},
 we can define the model set
$\mathcal{H}$ as
\begin{equation}\label{E:model_set}
\mathcal{H}\doteq\left\{ \theta\in\mathbb{R}^{n_\theta}:
\begin{array}{l}
0\leq f_{\hat{a}}(T,\theta_{\hat{a}}(\overline{p},\overline{y},\mu_i))<1,\,i=1,\ldots,M \\
0\leq
f_{\hat{e}}(T,\theta_{\hat{e}}(\overline{p},\overline{y},\mu_i))\leq
1,\,\forall
T\in\mathcal{T},\,i=1,\ldots,M\\
\sum\limits_{i=1}^M
f_{\hat{e}}(T,\theta_{\hat{e}}(\overline{p},\overline{y},\mu_i))=1\end{array},\;\;\forall
T\in\mathcal{T}\right\}
\end{equation}

As to the cost function $J$, this should represent a criterion by
which the accuracy of the ROM is evaluated. As discussed in section
\ref{S:SystPersp}, an indicator  of the accuracy of a ROM, suitable
for the considered application, is related to the error between its
total output intensity profile and that of the FOM, given the same
 temperature profile. Motivated by this consideration, we
select a series of temperature profiles, and then evaluate the
discrepancy between the corresponding outputs generated by the FOM
and those given by the ROM. More specifically, let us consider a
space interval $X=[0,\overline{x}]$, chosen such that
$\overline{x}=(L-1)\,\Delta x$ for some $L\in\mathbb{N}$. Similarly
to the choice of $\Delta x$, the value of $\overline{x}$ depends on
the problem at hand: a typical choice is three-four times the
largest space-scale of interest (see section \ref{SSS:input_design}
for details). Thus, in the interval $X$ we have $L$ position nodes
$x_l,\,l=1,\ldots,L$, with $x_1=0$ and $x_L=\overline{x}$. Consider
now a finite sequence of $L_T\in\mathbb{N}$ discretized temperature
profiles
$\tilde{T}_j\doteq[T_j(x_1),\ldots,T_j(x_L)]^T\in\mathbb{R}^{L\times
1},\,j=1,\ldots,L_T$. In face of each of such temperature profiles,
the FOM provides a corresponding profile of the radiated intensity,
$\tilde{I}_{tot,j}(\tilde{T}_j)\doteq[I_{tot,j}(x_1),\ldots,I_{tot,j}(x_L)]^T\in\mathbb{R}^{L\times
1},\,j=1,\ldots,L_T$, and similarly the ROM, for a given value of
$\theta$ provides an approximated one,
$\tilde{\hat{I}}_{tot,j}(\tilde{T}_j,\theta)$. Let us define the
weighted error profile as
\begin{equation}\label{E:cost_term}
\Delta I_{j}(\theta)\doteq
w_j(\tilde{I}_{tot,j}(\tilde{T}_j)-\tilde{\hat{I}}_{tot,j}(\tilde{T}_j,\theta)),\,j=1,\ldots,L_T,
\end{equation}
where $w_j>0$ is a scalar weight (the specific choice of $w_j$ is
discussed in section \ref{SSS:input_design}). Then, we take the cost
function as:
\begin{equation}\label{E:cost_func}
J(\theta)=\sum\limits_{j=1}^{L_T}\Delta I_{j}^T(\theta)\Delta
I_{j}(\theta),
\end{equation}
i.e. the sum, over all error profiles and over all position nodes of
each profile, of the weighted intensity error squared.

\subsubsection{Nonlinear Program formulation}\label{SS:NLP}
We can now write the model order reduction problem
 in a computationally tractable form:
\begin{equation}\label{E:pr_id_problem}
\min\limits_{\theta\in\mathcal{H}}J(\theta)
\end{equation}
where $\mathcal{H}$ and $J(\theta)$ are given in \eqref{E:model_set}
and \eqref{E:cost_func}, respectively. Namely, the aim of the
optimization problem \eqref{E:pr_id_problem} is to search, within
the model set $\mathcal{H}$, the value of the parameter vector
$\theta$ that achieves the best fitting criterion $J(\theta)$.
Problem \eqref{E:pr_id_problem} is a finite-dimensional nonlinear
program (NLP) which, for a suitable choice of the model
parametrization \eqref{E:parametrization}, can be tackled with
state-of-the-art numerical methods \cite{NoWr06}. As we show in the
examples of section \ref{S:results}, a typical dimension of the
optimization variable is $n_\theta\approx 77$, corresponding for
example to a 3$^\text{rd}$-order model with $25$ parameters defining
the functions $\hat{a}(T,\overline{p},\overline{y},\mu_i)$ and two
parameters defining the frequency cuts in the EM spectrum that
separate the three bands $\mu_1,\,\mu_2,\,\mu_3$. In general, the
cost function \eqref{E:cost_func} is non-convex with respect to the
decision variable $\theta$, so that only a local solution to problem
\eqref{E:pr_id_problem} can be computed efficiently with such
problem dimensions. As we show in the examples of section
\eqref{S:results}, the obtained solutions nevertheless
yield ROMs that are very accurate with respect to the FOM.\\
In the next section, we discuss several aspects involved in the
formulation and solution of \eqref{E:pr_id_problem}, including the
design of the input temperature profiles and the choice of the
weights in the cost function $J$ \eqref{E:cost_func}, the choice of
model parametrization, the numerical approach to solve the NLP,
finally some ways to improve the efficiency and stability of the
numerical optimization.

\subsection{Computational aspects}\label{SS:Computation}

\subsubsection{Input design, weighting coefficients and initial conditions}\label{SSS:input_design}
\begin{figure}[hbt]
\centerline{
\begin{tabular}{cc}
(a)&(b)\\
\includegraphics[clip,width=8cm]{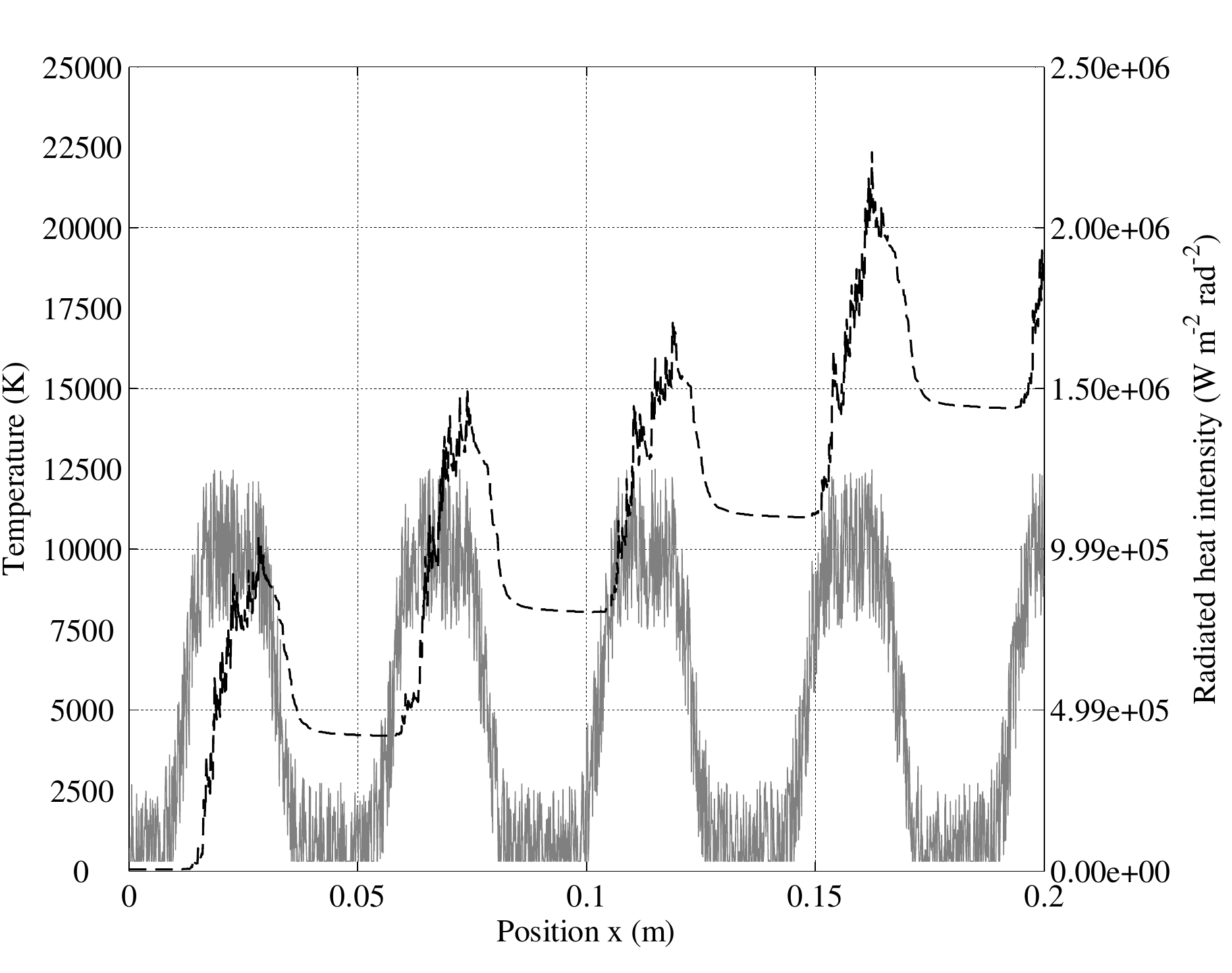}
&
\includegraphics[clip,width=8cm]{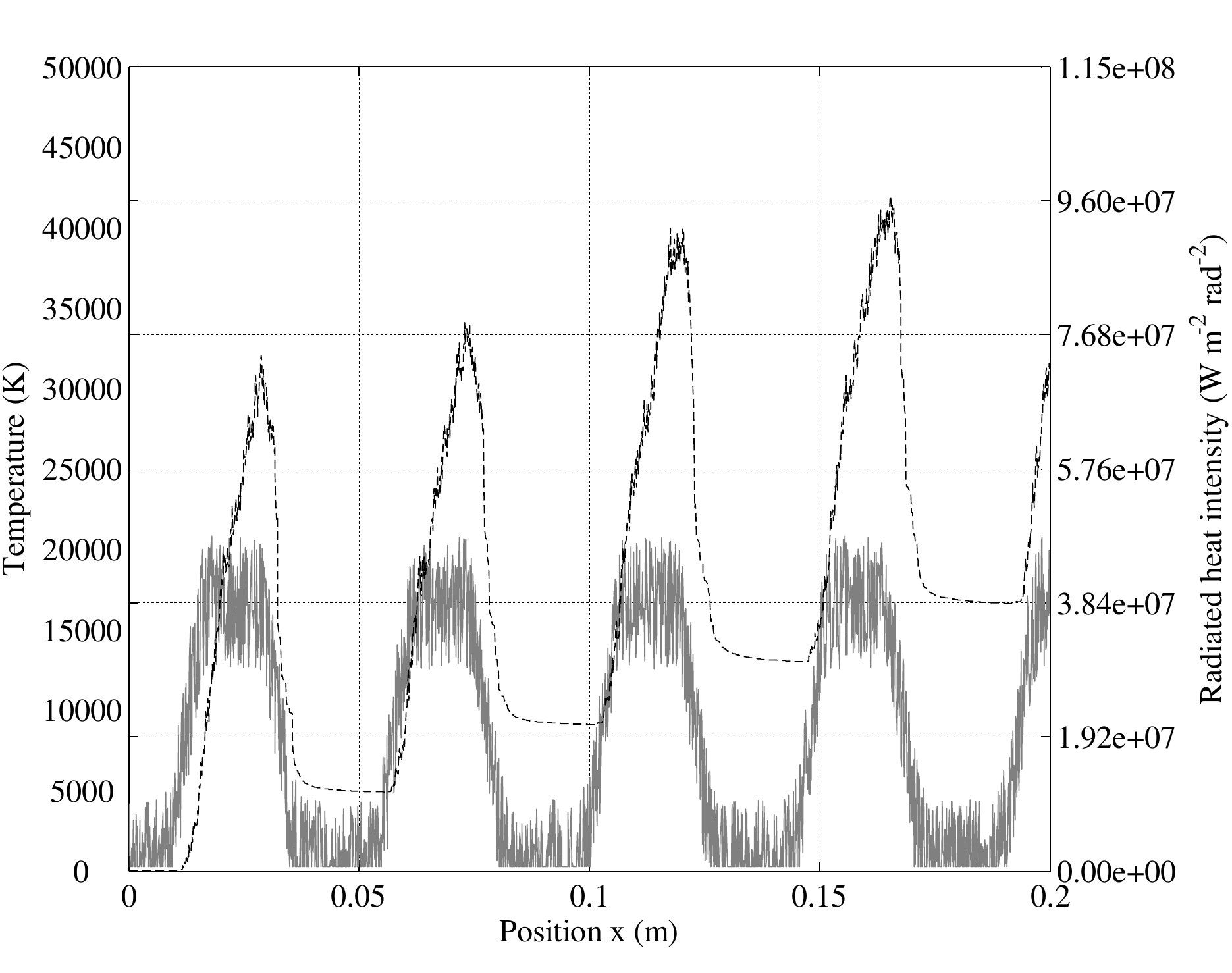}
\end{tabular}
} \protect\caption{\small Example of input temperature profiles used
in the order reduction procedure (solid gray lines) and radiated
heat intensity distribution computed with the full-order model
(dashed black lines). The maximum temperature values prior to adding
white noise are (a) $10^4\,$K  and (b) $1.7\;10^4\,$K . Base data:
pure air at $1.5\,10^5\,$Pa.
\normalsize}\label{F:bell_T_profile_noise}
\end{figure}

From an application's perspective, the position interval
$\overline{x}$, the resolution of the position discretization
$\Delta x$, and the temperature profiles
$\tilde{T}_j,\,j=1,\ldots,L_T$ should be representative of the
typical temperature distributions in the plasma that is generated
during the switching process. About the choice of $\Delta x$,
considering the discussion in section \ref{SS:space_discr}, a
reasonable tradeoff between accuracy and computational speed is
$\Delta x=10^{-4}\,$m, while for $\overline{x}$ a good choice in our
experience is $\overline{x}=2\,10^{-1}\,$m, for arcs whose width is
between $\approx3\,10^{-3}\,$m and $\approx2\,10^{-2}\,$m,
corresponding to current values in the range 1-40$\,10^3\,$A.
Finally, for the temperature profiles, we select a bell-shaped
function, where we can adjust the maximum temperature reached as
well as the steepness of the rising and falling slopes. Then, we
generate a series of such profiles by cycling through different
maximum temperatures and transient slopes. Furthermore, from a
system identification's perspective, given the nonlinearity of the
system's equations, different initial conditions for the internal
states of the model should be considered, and the input should be
designed in order to excite the system in a broad range of
frequencies. For the former aspect, we replicate the temperature
profile several times, so that the system is presented more than
once with the same temperature profile, but each time starting from
the internal states resulting from the previous profile. For the
second aspect, white noise processes are a well-known choice
\cite{Ljung99} to excite the system's dynamics over all frequencies,
hence we super-impose to the computed profiles a uniformly
distributed white noise temperature signal whose amplitude is a
fraction (e.g. $25\%$) of the highest temperature in the original
profile. With these choices, typical temperature profiles are shown
in Figure \ref{F:bell_T_profile_noise}(a)-(b), together with the
resulting intensity distributions computed with the FOM for a
mixture of pure air.

Due to the fourth-order dependence of the total black-body intensity
on temperature (see \eqref{E:total_bb}), the radiated intensity
obtained with different temperature profiles can have different
orders of magnitude, compare for example Figure
\ref{F:bell_T_profile_noise}(a) with Figure
\ref{F:bell_T_profile_noise}(b), where the maximum temperatures are
about 10$^4\,$K and $1.7\;10^4\,$K, respectively. If no corrective
measure is taken, such a difference among the various intensity
profiles would lead to poor accuracy of the reduced order models,
due to the biasing towards the higher intensity values in the cost
function. In order to compensate this effect, we select the weights
$w_j$ in \eqref{E:cost_term} on the basis of the radiated heat
intensity given by the FOM when the corresponding temperature
profile, $\tilde{T}_j$, is considered. In this way, all the error
profiles \eqref{E:cost_term} are normalized with respect to the
corresponding intensity; typical specific choices for the weights
include the maximum or the average radiated heat, i.e.:
\[
\begin{array}{c}
w_j=\max(\tilde{I}_{tot,j}),\,j=1,\ldots,L_T\\
\text{or}\\
w_j=\text{mean}(\tilde{I}_{tot,j}),\,j=1,\ldots,L_T\\
\end{array}
\]

Finally, it is worth mentioning the choice of the initial conditions
for the radiated heat intensity,  for both the FOM and the ROM, i.e.
the values of vectors $\overline{I}(x_1)$ and
$\hat{\overline{I}}(x_1)$ needed to simulate the models
\eqref{E:FOM_structure_dt} and \eqref{E:ROM_structure_dt},
respectively, in the interval $X=[0,\,\overline{x}]$. To this end,
we assume that for $x<0$ (i.e. outside such an interval) the
incoming radiation in positive $x$ direction corresponds to the
steady-state intensity at ambient temperature $T_a$, i.e.
$I_{bb}(T_a,\nu_i),\,i=1,\ldots,N$ for the FOM and
\begin{equation}\label{E:init_cond_ROM}
\hat{I}_{j}(x,\mu_i)=\hat{e}(T_a,\overline{p},\overline{y},\mu_i)\overline{I}_{bb}(T_a),\,i=1,\ldots,M
\end{equation}
for the ROM. The rationale behind this choice is that in the arc
simulations related to switchgear devices one can assume that
initially the boundaries of the computational domain are at ambient
temperature and radiate the corresponding intensity into the volume
filled with plasma. Additionally, the boundaries can not reach
temperatures higher than about 3,000 K (depending on the material),
such that the related radiated heat is anyway negligible with
respect to the intensity emitted in the plasma volume where current
is flowing.

\subsubsection{Model parametrization and model set}\label{SSS:parametrization}
The choice of the model parametrization, i.e. of functions
$f_{\hat{a}}$ and $f_{\hat{e}}$ in \eqref{E:parametrization}, is
crucially important for the accuracy of the obtained results and for
the solution of the numerical optimization. For the first aspect,
one shall choose a rich enough family of functions, such that the
data from the FOM can be reproduced with small errors. For the
second aspect, the best choice would be a parametrization leading to
a convex model set \eqref{E:model_set}, so that sequential quadratic
programming and line search algorithms \cite{NoWr06} can be used
efficiently. For the equivalent absorption coefficients $\hat{a}$,
i.e. function $f_{\hat{a}}$, a choice that meets both requirements
is a piecewise affine parametrization, where
$\theta_{\hat{a}}(\overline{p},\overline{y},\mu_i)$ is, for each
band $i=1,\ldots,M$, a vector containing the values of the
coefficients $\hat{a}$ at a finite number of pre-defined temperature
nodes $T_k,\,k=1,\ldots,N_k$, such that
$T_{k-1}<T_k,\,k=2,\ldots,N_k$, chosen by the user (e.g. equally
spaced):
\[
\theta_{\hat{a}}(\overline{p},\overline{y},\mu_i)=\left[\begin{array}{c}
\theta_{\hat{a}}(T_1,\overline{p},\overline{y},\mu_i)\\
\vdots\\
\theta_{\hat{a}}(T_{N_k},\overline{p},\overline{y},\mu_i)
\end{array}\right]\in\mathbb{R}^{N_k},\;\;i=1,\ldots,M.
\]
The temperature nodes must include the values at the boundaries
\eqref{E:sets} of the domain $\mathcal{T}$, i.e. $T_1=T_\text{min}$
and
$T_{N_k}=T_\text{max}$.\\
Then, for a given temperature value $T\in\mathcal{T}$, the function
$f_{\hat{a}}$ in \eqref{E:parametrization} is computed by
interpolating linearly among the values of $\theta_{\hat{a}}$
corresponding to the neighboring temperature nodes:
\begin{equation}\label{E:func_param_pwaffine}
f_{\hat{a}}(T,\theta_{\hat{a}}(\overline{p},\overline{y},\mu_i))=\lambda(T)^T\theta_{\hat{a}}(\overline{p},\overline{y},\mu_i),
\end{equation}
where \begin{equation}\label{E:lambda} \lambda(T)=
\left[\begin{array}{c}
0\\\vdots\\1-\dfrac{T-T_{k-1}}{T_k-T_{k-1}}\\\dfrac{T-T_{k-1}}{T_k-T_{k-1}}\\\vdots\\0
\end{array}\right]\in\mathbb{R}^{N_k}
\end{equation}
and $T_{k-1},\,T_{k}$ are two subsequent nodes such that
$T\in[T_{k-1},\,T_{k}]$. Using the parametrization
\eqref{E:func_param_pwaffine}, the number of optimization variables
introduced in the problem is $M\,N_k$. The flexibility of the
approximating function $f_{\hat{a}}$  can be increased by increasing
the number of temperature nodes $N_k$. Moreover, note that
$f_{\hat{a}}$ in \eqref{E:func_param_pwaffine} is a convex
combination \cite{BoVa04} of the elements contained in the parameter
vector. Therefore, with this parametrization,  the constraints
\eqref{E:constraints_dt_func} can be enforced by imposing them just
on the values at the nodes, i.e.:
\begin{equation}\label{E:constraints_parameters_a}
\begin{array}{l}
0\leq
\hat{a}(T_k,\overline{p},\overline{y},\mu_i)<1,\,k=1,\ldots,N_k,\,i=1,\ldots,M
\end{array}.
\end{equation}

As regards the parametrization of the equivalent emissivity
functions, $\hat{e}(T,\overline{p},\overline{y},\mu_i)$, an
additional requirement is to retain a physical link between the
emissivities of each band in the ROM and the original absorption
spectrum, such that each $\mu_i$ accounts for the intensity
contributed by a precise frequency interval defined by suitable
frequency cuts. As already mentioned, previous contributions in the
literature actually consider such a partitioning of the EM spectrum
\cite{NoIo08,JCGB14}, however the choice of the frequency bands is
not trivial and not systematic, so that one has to proceed with a
trial-and-error approach. With the technique proposed here, one can
optimize directly with respect to such frequency cuts, hence
obtaining a systematic method for the band-averaging. This can be
done by using the following model parametrization for the equivalent
emissivities:
\begin{equation}\label{E:func_param_emis_bands}
f_{\hat{e}}(T,\theta_{\hat{e}}(\overline{p},\overline{y},\mu_i))=\gamma(T,\theta_{\hat{e}}(\overline{p},\overline{y},\mu_i))\doteq\int\limits_{\theta_{\hat{e}_{hf}}
(\overline{p},\overline{y},\mu_{i-1})}^{\theta_{\hat{e}_{hf}}(\overline{p},\overline{y},\mu_i)}e(T,\nu)d\nu,
\end{equation}
where $e(T,\nu)$ is the spectral emissivity \eqref{E:spec_emis} and
the optimization variables are the upper boundaries (e.g. in Hz) of
each frequency band:
\[
\theta_{\hat{e}}(\overline{p},\overline{y},\mu_i)=\left[\begin{array}{l}
\theta_{\hat{e}_{hf}}(\overline{p},\overline{y},\mu_1)\\\vdots\\
\theta_{\hat{e}_{hf}}(\overline{p},\overline{y},\mu_{M-1})
\end{array}\right]\in\mathbb{R}^{M-1}.
\]
For the first and last bands, i.e. $\mu_1$ and $\mu_M$, the
equivalent emissivity is computed using
\eqref{E:func_param_emis_bands} with $\mu_0\rightarrow0$ and
$\mu_M=+\inf$, respectively.
 In this way, the equivalent emissivity of
each band is, by construction, equal to the one pertaining to the
frequency interval
$[\theta_{\hat{e}_{hf}}(\overline{p},\overline{y},\mu_{i-1}),\,\theta_{\hat{e}_{hf}}(\overline{p},\overline{y},\mu_i)]$.
The constraints \eqref{E:constraints_ct_func_e} defining the model
set (i.e. the second and third rows in \eqref{E:model_set}) can be
easily taken into account by a set of linear inequalities:
\begin{equation}\label{E:constraints_parameters_e_bands}
\begin{array}{l}
\theta_{hf}(\overline{p},\overline{y},\mu_{1})>0\\
\theta_{hf}(\overline{p},\overline{y},\mu_{i+1})>\theta_{hf}(\overline{p},\overline{y},\mu_{i}),\,i=2,\ldots,M-2
\end{array}.
\end{equation}
The constraints \eqref{E:constraints_parameters_e_bands} impose an
increasing ordering of the boundaries. By construction, the
resulting bands
 are not overlapping and they cover the
whole spectrum: these features automatically enforce  the
constraints \eqref{E:constraints_ct_func_e}. We note that with this
parametrization, the number of free variables that define the band
partition is equal to $M-1$. Therefore, considering also the
parameters pertaining to the equivalent absorption coefficients, the
total number of optimization variables in the problem is equal to
$M\,N_k+M-1$. About the total number of constraints, in virtue of
\eqref{E:constraints_parameters_a} we have a set of $2\,M\,N_k$
linear inequalities, while \eqref{E:constraints_parameters_e_bands}
amounts to $M-1$ additional inequalities, for a total of
$2\,M\,N_k+M-1$. These inequalities altogether define indeed a
convex set (in particular a polytope) where the solution of the
optimization problem is confined.

About the computation of
$\gamma(T,\theta_{\hat{e}}(\overline{p},\overline{y},\mu_i))$ in
\eqref{E:func_param_emis_bands}, note that this function can be also
written as:
\begin{equation}\label{E:f_gamma_def}
\gamma(T,\theta_{\hat{e}}(\overline{p},\overline{y},\mu_i))=f_\gamma(T,\theta_{\hat{e}_{hf}}
(\overline{p},\overline{y},\mu_i))-f_\gamma(T,\theta_{\hat{e}_{hf}}(\overline{p},\overline{y},\mu_{i-1}))
\end{equation}
where
\begin{equation}\label{E:f_gamma}
f_\gamma(T,\theta)\doteq\int\limits_0^{\theta}e(T,\nu)d\nu.
\end{equation}
The function
$f_\gamma:\mathcal{T}\times\mathbb{R}^+\rightarrow[0,1]$ can be
conveniently pre-computed and stored, so that the computation of
$\gamma(T,\theta_{\hat{e}}(\overline{p},\overline{y},\mu_i))$ can be
very efficient.

\subsubsection{NLP solution and computational aspects}\label{SSS:comp_eff}
After choosing the temperature profiles and the model
parametrization as described in sections \ref{SSS:input_design} and
\ref{SSS:parametrization}, respectively, it can be shown that the
cost function $J$ \eqref{E:cost_func} is twice differentiable. Thus,
considering also the convexity of the model set, the optimization
problem \eqref{E:pr_id_problem} can be solved for a local minimum
with constrained sequential quadratic programming (SQP) techniques
\cite{NoWr06}, like the one implemented in Matlab$^\circledR$
function \verb"fmincon". The computational efficiency and stability
of the numerical optimization depend on a number of aspects, which
are briefly mentioned here.

\emph{Computation of the cost function}. SQP solvers employ (as most
NLP solution algorithms) an iterative strategy where the cost
function might need to be evaluated hundreds of times. It is thus
imperative to speed-up the computation of such a function which, in
our case, implies the computation of the $L-$dimensional vectors
$\Delta I_{j}(\theta)=$
$w_j(\tilde{I}_{tot,j}(\tilde{T}_j)-\tilde{\hat{I}}_{tot,j}(\tilde{T}_j,\theta)),\,j=1,\ldots,L_T$,
see \eqref{E:cost_term}-\eqref{E:cost_func}. The intensity profiles
given by the FOM, $\tilde{I}_{tot,j}(\tilde{T}_j)$, can be computed
once and used in all the function evaluations, since they are not
changing with $\theta$. On the other hand, the intensity given by
the ROM, $\tilde{\hat{I}}_{tot,j}(\tilde{T}_j,\theta)$, has to be
computed iteratively during the optimization. This implies
evaluating the selected functions \eqref{E:parametrization} and
simulating the ROM model. As mentioned briefly in section
\ref{SS:space_discr}, the computational times of such operations
improve significantly by using the discretized version of the ROM
\eqref{E:ROM_structure_dt} and by considering the values of
$\hat{a}$, instead of $\hat{\alpha}$, as optimization parameters,
since one avoids the computation of the exponential in
\eqref{E:coeff_a_dt}. The values of $\hat{\alpha}$ can be then
recovered from the optimal solution by inverting such equation.
Another approach to speed up the computations is to pre-compute and
store the vectors
 $\lambda(T_j(x_i)),\,i=1,\ldots,L,\,j=1,\ldots,L_T$ \eqref{E:lambda} corresponding to the employed
temperature profiles, which are needed to compute the piecewise
affine functions
$f_{\hat{a}}(T,\theta_{\hat{a}}(\overline{p},\overline{y},\mu_i))$
\eqref{E:func_param_pwaffine}, since these vectors do not depend on
the optimization variables but just on the temperature distributions
$\tilde{T}_j$.

\emph{Gradient computations}. Another aspect that greatly influences
the efficiency (and accuracy) of the optimization algorithm is the
computation of the gradients of the cost and constraint functions.
While the gradient of each constraint is trivial to compute, since
only linear equalities and inequalities are present (as shown in
section \ref{SSS:parametrization}), the gradient of the cost is more
difficult to obtain, but we show here that it can still be derived
analytically. First of all, consider that:
\begin{equation}\label{E:gradient1}
J(\theta)=\sum\limits_{j=1}^{L_T}\sum\limits_{l=1}^L  w_j^2
(I_{tot,j}(x_l)-\hat{I}_{tot,j}(x_l,\theta))^2
\end{equation}
hence
\begin{equation}\label{E:gradient2}
\nabla_\theta J(\theta)\doteq\left[
\begin{array}{c}
\frac{\partial J}{\partial
\theta_{\hat{a},1}(\overline{p},\overline{y},\mu_1)}\\
\vdots\\
\frac{\partial J}{\partial
\theta_{\hat{e},n_{\theta_{\hat{e}}}}(\overline{p},\overline{y},\mu_M)}
\end{array}\right]=\sum\limits_{j=1}^{L_T}\sum\limits_{l=1}^L
-2\,w_j^2\left(I_{tot,j}(x_l)-\hat{I}_{tot,j}(x_l,\theta)\right)\nabla_\theta\hat{I}_{tot,j}(x_l,\theta).
\end{equation}
We thus focus on deriving a general expression for the gradient
$\nabla_\theta\hat{I}_{tot,j}(x_l,\theta)$ of the total radiated
intensity at the $l-$th position node for the $j-$th temperature
profile, which can then be employed to compute the gradient of $J$
using \eqref{E:gradient2}. Exploiting the second of the model
equations \eqref{E:ROM_structure_dt} and considering the vector
$\hat{C}$ in \eqref{E:ROM_matrices}, we have:
\begin{equation}\label{E:gradient3}
\nabla_\theta
\hat{I}_{tot,j}(x_l,\theta)=\sum\limits_{i=1}^M\nabla_\theta\hat{I}_{j}(x_l,\mu_i,\theta),
\end{equation}
where $\hat{I}_{j}(x_l,\mu_i,\theta)$ is the intensity pertaining to
the band $\mu_i$ at the position step $x_l$ when the temperature
profile $\tilde{T}_j$ is considered. Moreover, by using the first of
\eqref{E:ROM_structure_dt} together with \eqref{E:ROM_matrices_dt}
and \eqref{E:parametrization} we obtain, for each $l=2,\ldots,L$:
\begin{equation}\label{E:gradient4}
\begin{array}{rcl}
\nabla_\theta\hat{I}_{j}(x_l,\mu_i,\theta)&=&\nabla_\theta
f_{\hat{a}}(T_j(x_{l-1}),\theta_{\hat{a}}(\overline{p},\overline{y},\mu_i))\;\hat{I}_{j}(x_{l-1},\mu_i,\theta)\\
&+&f_{\hat{a}}(T_j(x_{l-1}),\theta_{\hat{a}}(\overline{p},\overline{y},\mu_i))\;\nabla_\theta\hat{I}_{j}(x_{l-1},\mu_i,\theta)\\
&+&\left(1-f_{\hat{a}}(T_j(x_{l-1}),\theta_{\hat{a}}(\overline{p},\overline{y},\mu_i))\right)\;\nabla_\theta
f_{\hat{e}}(T_j(x_{l-1}),\theta_{\hat{e}}(\overline{p},\overline{y},\mu_i))\;\overline{I}_{bb}(T_j(x_{l-1}))\\
&-&\nabla_\theta
f_{\hat{a}}(T_j(x_{l-1}),\theta_{\hat{a}}(\overline{p},\overline{y},\mu_i))\;
f_{\hat{e}}(T_j(x_{l-1}),\theta_{\hat{e}}(\overline{p},\overline{y},\mu_i))\;\overline{I}_{bb}(T_j(x_{l-1}))
\end{array}
\end{equation}
The gradients $\nabla_\theta f_{\hat{a}}$ and $\nabla_\theta
f_{\hat{e}}$ in \eqref{E:gradient4} depend on the model
parametrization. For the piecewise affine function $f_{\hat{a}}$,
the computation is straightforward, e.g. (from
\eqref{E:func_param_pwaffine}):
\begin{equation}\label{E:gradient5}
\nabla_\theta
f_{\hat{a}}(T_j(x_{l-1}),\theta_{\hat{a}}(\overline{p},\overline{y},\mu_i)))=\nabla_\theta\left(\lambda(T_j(x_{l-1}))\theta_{\hat{a}}(\overline{p},\overline{y},\mu_i)\right)=
\left[\begin{array}{c} 0\\
\vdots\\
\lambda(T_j(x_{l-1}))\\
\vdots\\0
\end{array}\right]\in\mathbb{R}^{n_\theta},
\end{equation}
where the zeros in vector $\nabla_\theta f_{\hat{a}}$ correspond to
all the other parameters in vector $\theta$ but the set
$\theta_{\hat{a}}(\overline{p},\overline{y},\mu_i)$ pertaining to
the $i^\text{th}$ band.\\
As regards function $f_{\hat{e}}$, whose parameters are the
boundaries of the spectral bands, the computation of the related
gradient $\nabla_\theta f_{\hat{e}}$ involves computing the
derivative of function $f_\gamma(T,\theta)$ in \eqref{E:f_gamma},
which is readily obtained as
$df_\gamma(T,\theta)/d\theta=e(T,\theta)$ \eqref{E:spec_emis}. More
specifically, considering \eqref{E:func_param_emis_bands} and
\eqref{E:f_gamma_def}-\eqref{E:f_gamma} we have:
\begin{equation}\label{E:gradient6}
\begin{array}{rcl}
\nabla_\theta
f_{\hat{e}}(T_j(x_{l-1}),\theta_{\hat{e}}(\overline{p},\overline{y},\mu_i)))&=&
\nabla_\theta\gamma\left(T_j(x_{l-1}),\theta_{\hat{e}}(\overline{p},\overline{y},\mu_i)\right)\\
&=&
\left[\begin{array}{c} 0\\
\vdots\\
-e(T_j(x_{l-1}),\theta_{\hat{e}_{hf}}(\overline{p},\overline{y},\mu_{i-1}))\\
e(T_j(x_{l-1}),\theta_{\hat{e}_{hf}}(\overline{p},\overline{y},\mu_{i}))\\
\vdots\\0
\end{array}\right]\in\mathbb{R}^{n_\theta},
\end{array}
\end{equation}
where, in a way similar to \eqref{E:gradient5}, the zeros in vector
$\nabla_\theta f_{\hat{e}}$ correspond to all the other parameters
in vector $\theta$ but the pair
$\theta_{\hat{e}_{hf}}(\overline{p},\overline{y},\mu_{i-1}),\,\theta_{\hat{e}_{hf}}(\overline{p},\overline{y},\mu_i)$
pertaining to the $i^\text{th}$ band. For the first and last bands,
i.e. $i=1$ and $i=M$ respectively, the gradient is computed by
replacing the spectral emissivity pertaining to
$\theta_{\hat{e}_{hf}}(\overline{p},\overline{y},\mu_{i-1})$ (resp.
$\theta_{\hat{e}_{hf}}(\overline{p},\overline{y},\mu_{i})$) with
zero, since these two boundaries are not free variables, as
commented above in Section \ref{SSS:parametrization}.

 The last ingredient needed for the gradient
computation is the initialization of
$\nabla_\theta\hat{I}_{j}(x_l,\mu_i,\theta)$ for $l=1$, which is
required in \eqref{E:gradient4}; this is readily done by computing
the derivative of the initial conditions \eqref{E:init_cond_ROM}:
\begin{equation}\label{E:gradient7}
\nabla_\theta\hat{I}_{j}(x_1,\mu_i,\theta)=\nabla_\theta
f_{\hat{e}}(T_a,\theta_{\hat{e}}(\overline{p},\overline{y},\mu_i))\;\overline{I}_{bb}(T_a),\,i=1,\ldots,M,\,j=1,\ldots,L_T
\end{equation}
By using \eqref{E:gradient4}-\eqref{E:gradient7}, one can
recursively compute the values of
$\nabla_\theta\hat{I}_{j}(x_l,\mu_i,\theta)$ for each position node
$x_l$, each temperature profile $\tilde{T}_j$ and each band $\mu_i$,
and then use \eqref{E:gradient1}-\eqref{E:gradient3} to compute the
gradient of the cost function. The recursive gradient computation
can be done together with the simulation  of the ROM, which is
anyway needed to compute the cost function. The alternative to such
an approach would be to estimate the gradient by means of numerical
differentiation, e.g. finite difference approximation, which would
imply computing $n_\theta+1$ times the cost function for a single
gradient estimate. This is done by default in most available
optimization routines if no gradient is provided, however it is
subject to numerical errors and it requires significantly more time
than the exact computation derived above. Therefore, the use of
\eqref{E:gradient1}-\eqref{E:gradient7} to compute $\nabla_\theta
J(\theta)$ is both more accurate and more efficient than numerical
differentiation.
\begin{figure}[hbt]
\centerline{
\includegraphics[clip,width=14cm]{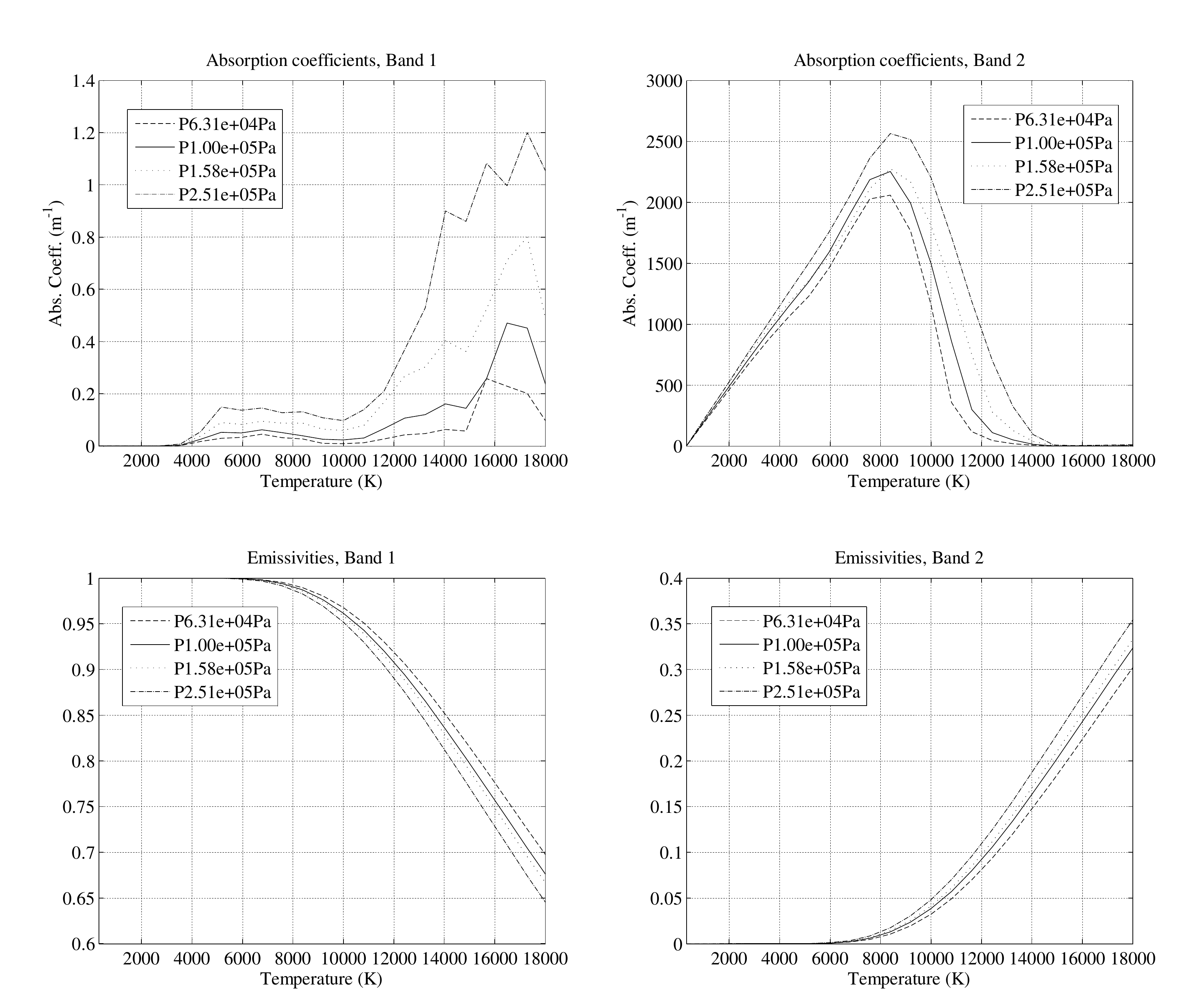}
} \protect\caption{\small Example of absorption coefficients and
emissivities for a ROM with two bands for pure
air.\normalsize}\label{F:abs_example1}
\end{figure}
\begin{figure}[hbt]
\centerline{
\includegraphics[clip,width=14cm]{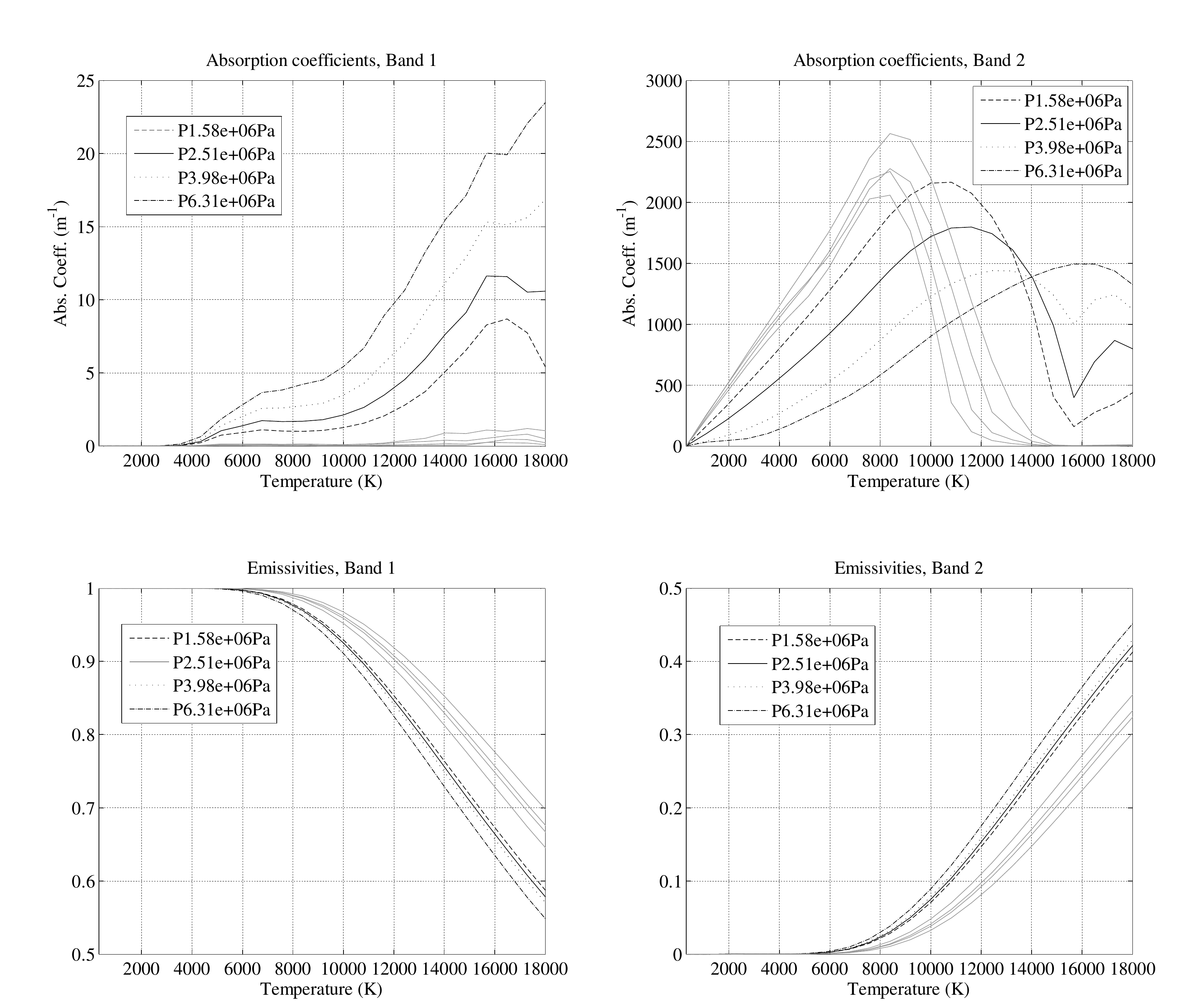}
} \protect\caption{\small Example of absorption coefficients and
emissivities for a ROM with two bands for pure air. The gray lines
correspond to the curves shown in Figure
\ref{F:abs_example1}\normalsize}\label{F:abs_example2}
\end{figure}

\begin{figure}[hbt]
\centerline{
\begin{tabular}{cc}
(a)&(b)\\
\includegraphics[clip,width=7.5cm]{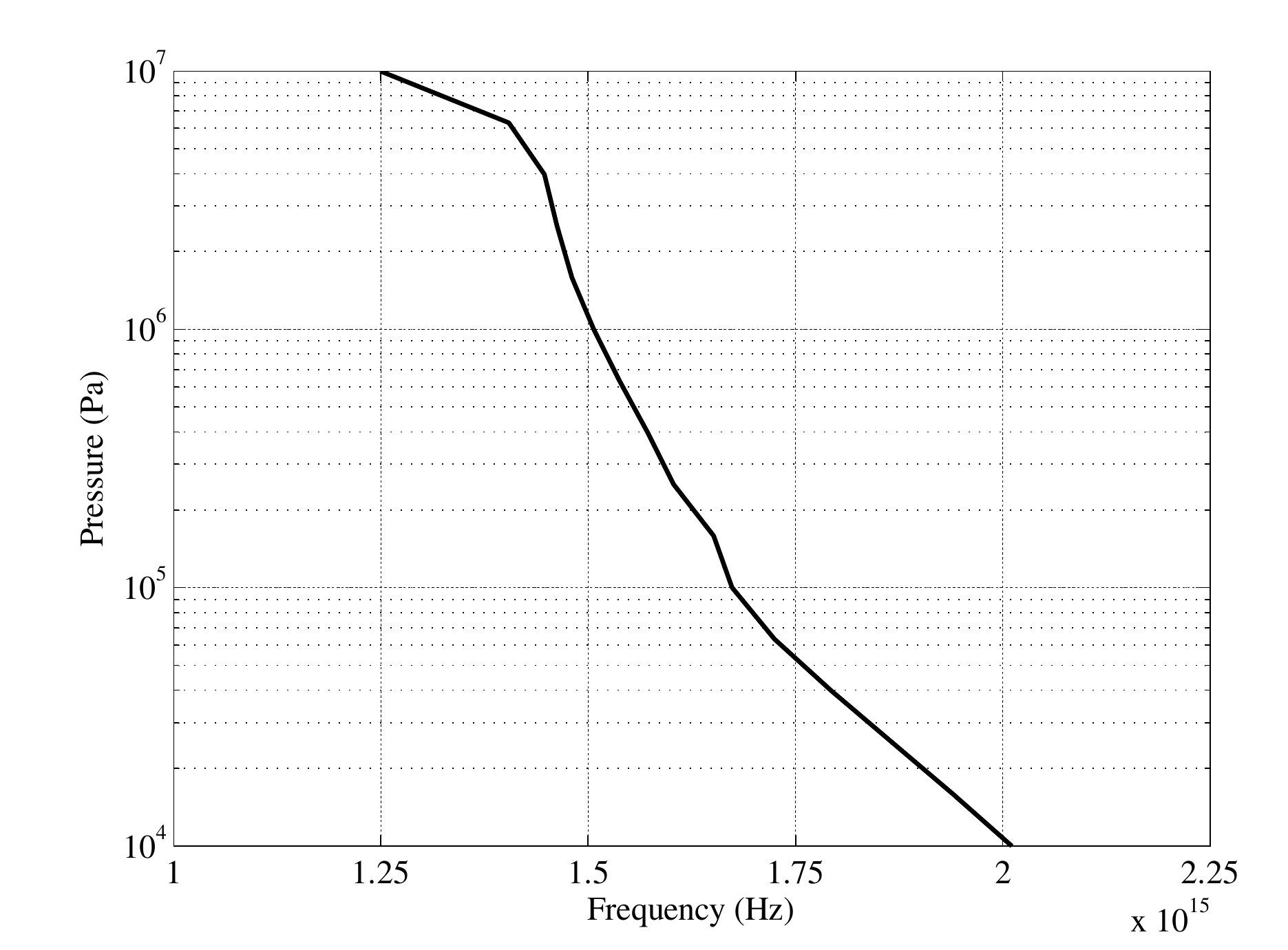}
&
\includegraphics[clip,width=7.5cm]{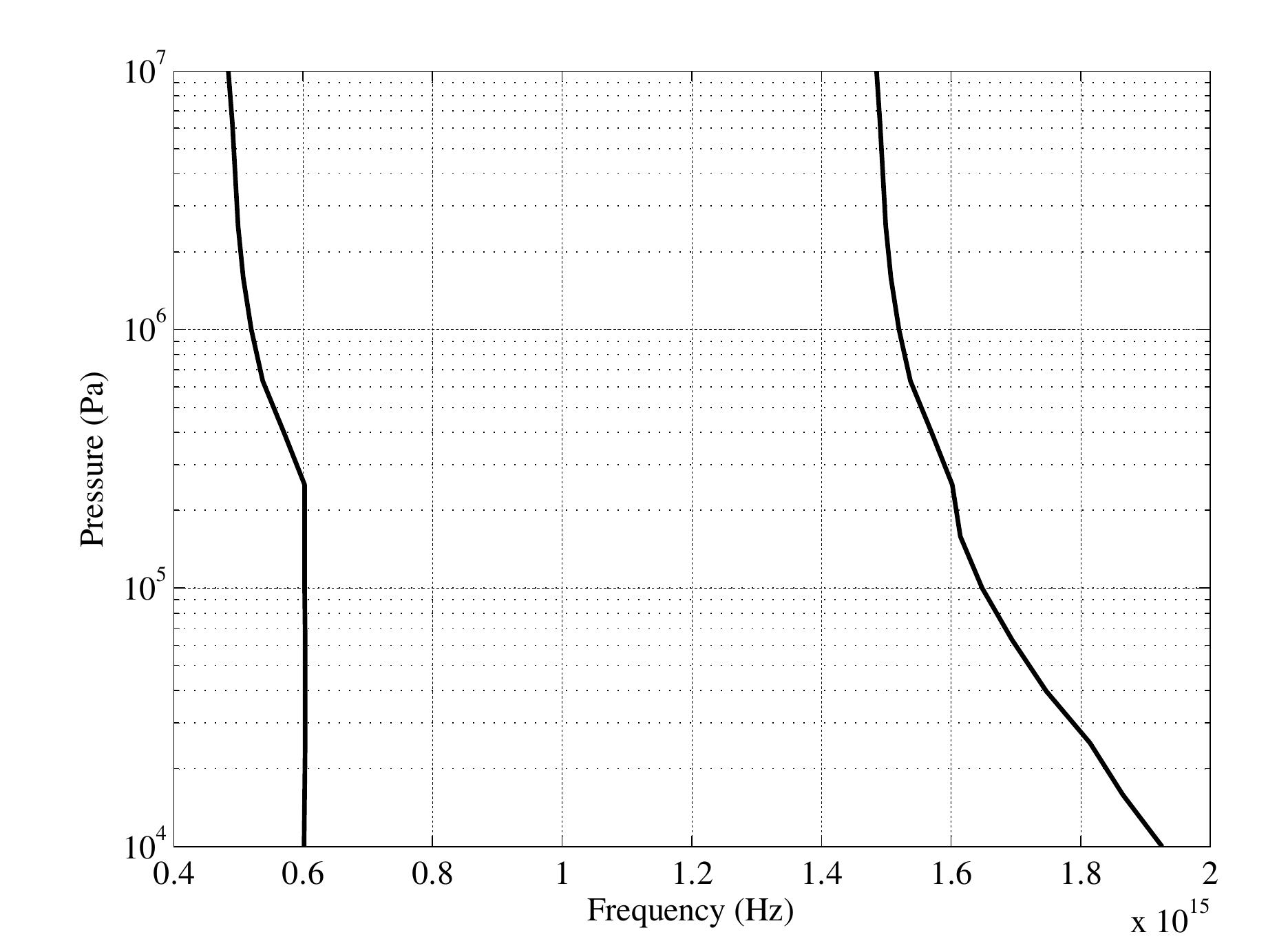}
\end{tabular}
} \protect\caption{\small Optimal partitioning of the frequency
spectrum for pure air and (a) two bands, (b) three
bands.\normalsize}\label{F:Opt_bands}
\end{figure}

\begin{figure}[hbt]
\centerline{
\includegraphics[clip,width=15cm]{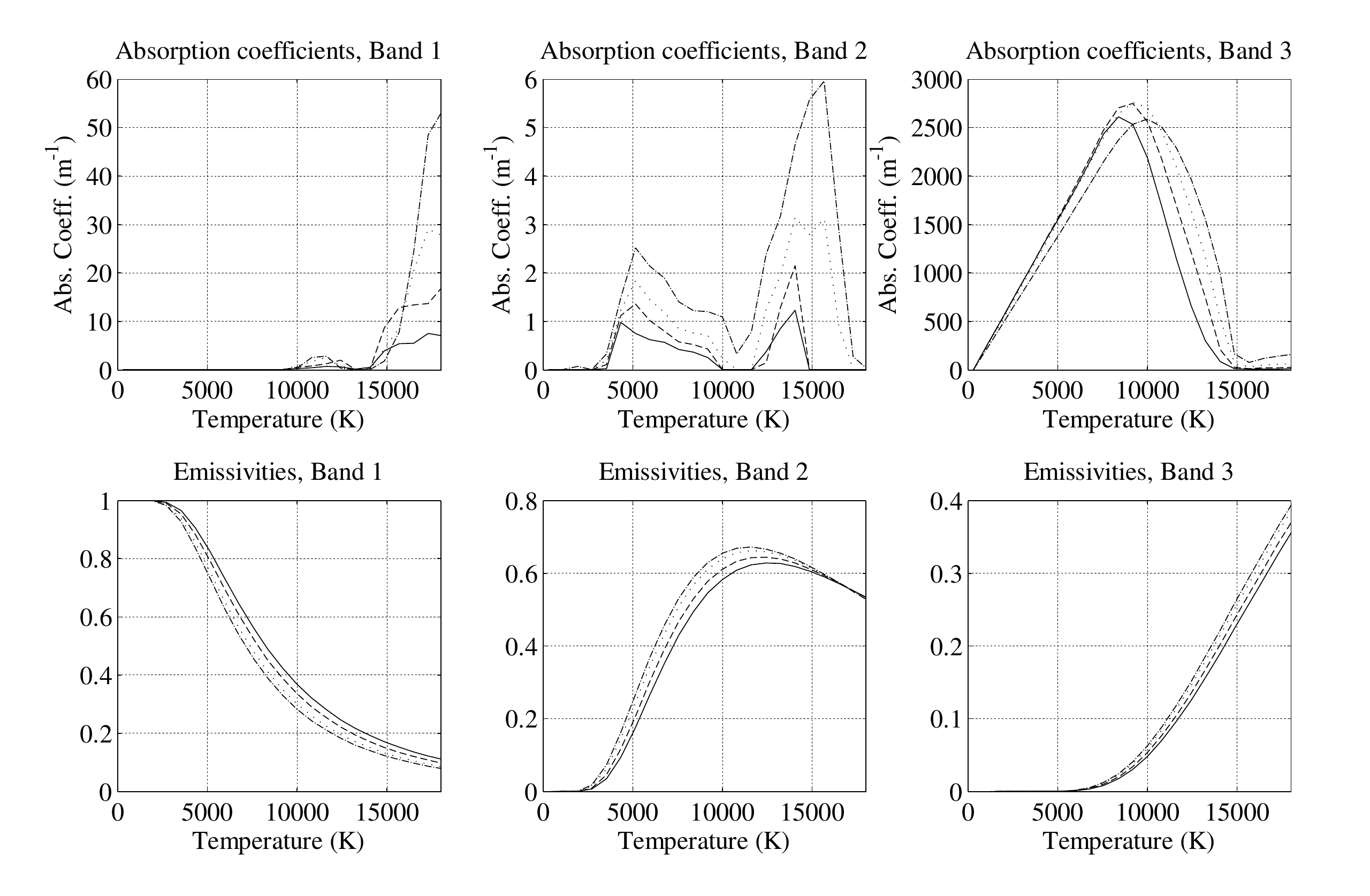}
} \protect\caption{\small Example of absorption coefficients and
emissivities for a ROM with three bands for pure air. Pressure
values: $2.51\,10^5\,$Pa (solid lines), $3.98\,10^5\,$Pa (dashed
lines), $6.31\,10^5\,$Pa (dotted lines), $1\,10^6\,$Pa (dash-dotted
lines). \normalsize}\label{F:abs_example3}
\end{figure}

\emph{Regularization and warm start}. Other important issues that
can arise in the numerical solution of \eqref{E:pr_id_problem} are
related to the quality of the obtained optimizer and the stability
of the optimization algorithm. One aspect that is relevant for the
subsequent use of the derived absorption coefficients in CFD
simulations is the smoothness of the obtained functions
$\hat{\alpha}$ and $\hat{\theta}$ with respect to their input
arguments, $T,p,y,\mu_i$. About the temperature dependence, if
piecewise affine functions are used one can penalize the variation
of $\hat{\alpha}$ (or equivalently $\hat{a}$) directly in the
optimization, by augmenting the cost function with a regularization
term, as follows:
\[
J(\theta)=\sum\limits_{j=1}^{L_T}\Delta I_{j}^T(\theta)\Delta
I_{j}(\theta)+\beta\Delta\theta^T\Delta\theta,
\]
where
\[
\Delta\theta\doteq\left[
\begin{array}{c}
\theta_{\hat{a},2}(\overline{p},\overline{y},\mu_1)-\theta_{\hat{a},1}(\overline{p},\overline{y},\mu_1)\\
\vdots\\
\theta_{\hat{a},n_{\theta_{\hat{a}}}}(\overline{p},\overline{y},\mu_1)-\theta_{\hat{a},n_{\theta_{\hat{a}}}-1}(\overline{p},\overline{y},\mu_1)\\
\vdots\\
\theta_{\hat{a},2}(\overline{p},\overline{y},\mu_M)-\theta_{\hat{a},1}(\overline{p},\overline{y},\mu_M)\\
\vdots\\
\theta_{\hat{a},n_{\theta_{\hat{a}}}}(\overline{p},\overline{y},\mu_M)-\theta_{\hat{a},n_{\theta_{\hat{a}}}-1}(\overline{p},\overline{y},\mu_M)
\end{array}\right]
\]
and $\beta$ is a weighting factor that can be used to tradeoff the
variability of the absorption coefficients with respect to
temperature with the goodness of the fitting criterion. The gradient
of this modified cost can be computed with a straightforward
extension of equation \eqref{E:gradient2}. About the smoothness with
respect to $p$ and $y$, recall that we assumed that several ROMs are
computed independently by gridding the $\cal{P}$ and $\cal{Y}$
domains (see section \ref{SS:Solution}). Thus, a possible way to
influence the dependence of the coefficients on pressure and
composition is to warm-start the optimization, by initializing the
optimization parameters with the solution obtained from neighboring
values of $\overline{p},\overline{y}$. This approach usually yields
good results, as we show in the next section with some examples.
\begin{figure}[hbt]
\centerline{
\begin{tabular}{cc}
(a)&(b)\\
\includegraphics[clip,width=7cm]{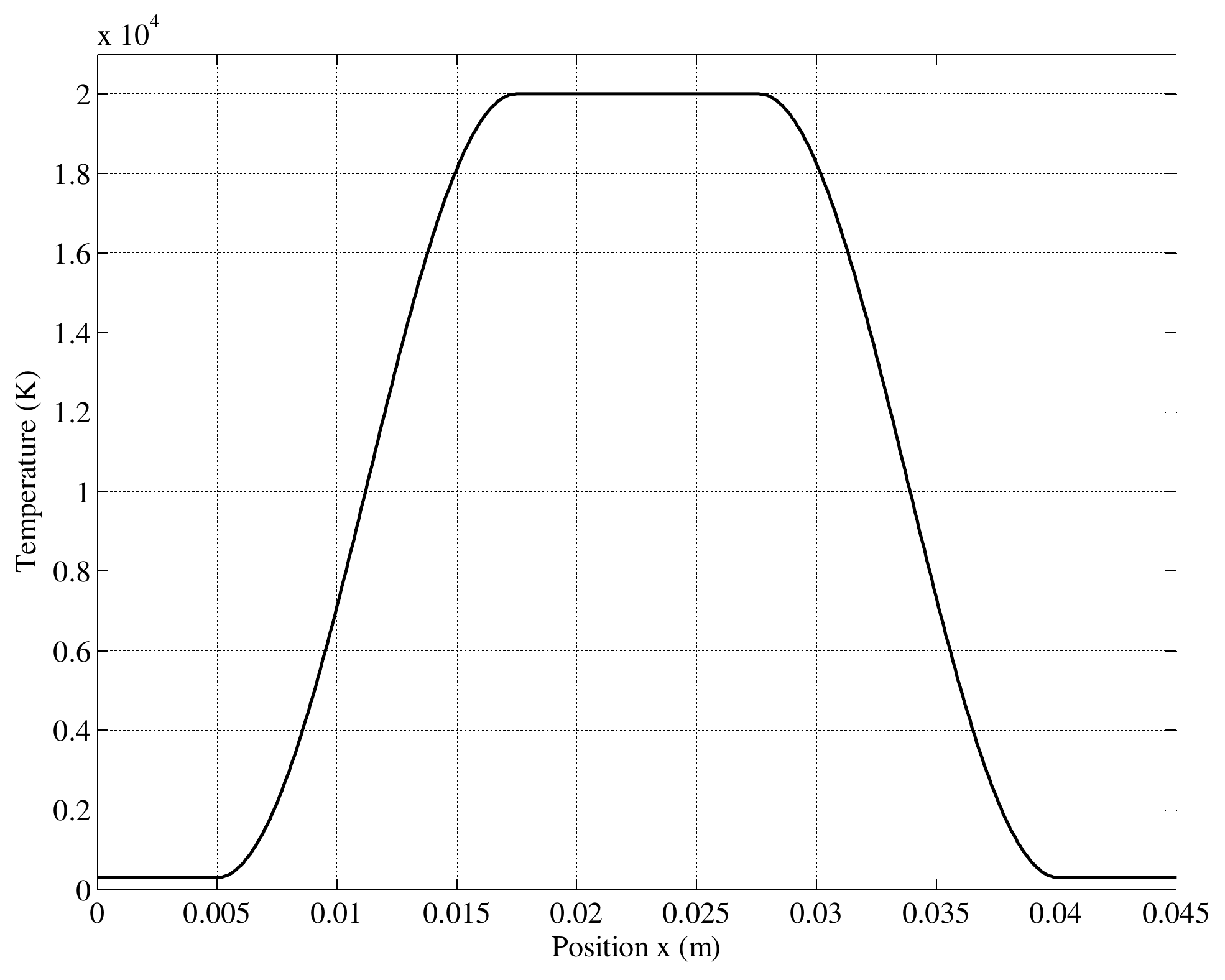}
&
\includegraphics[clip,width=7.5cm]{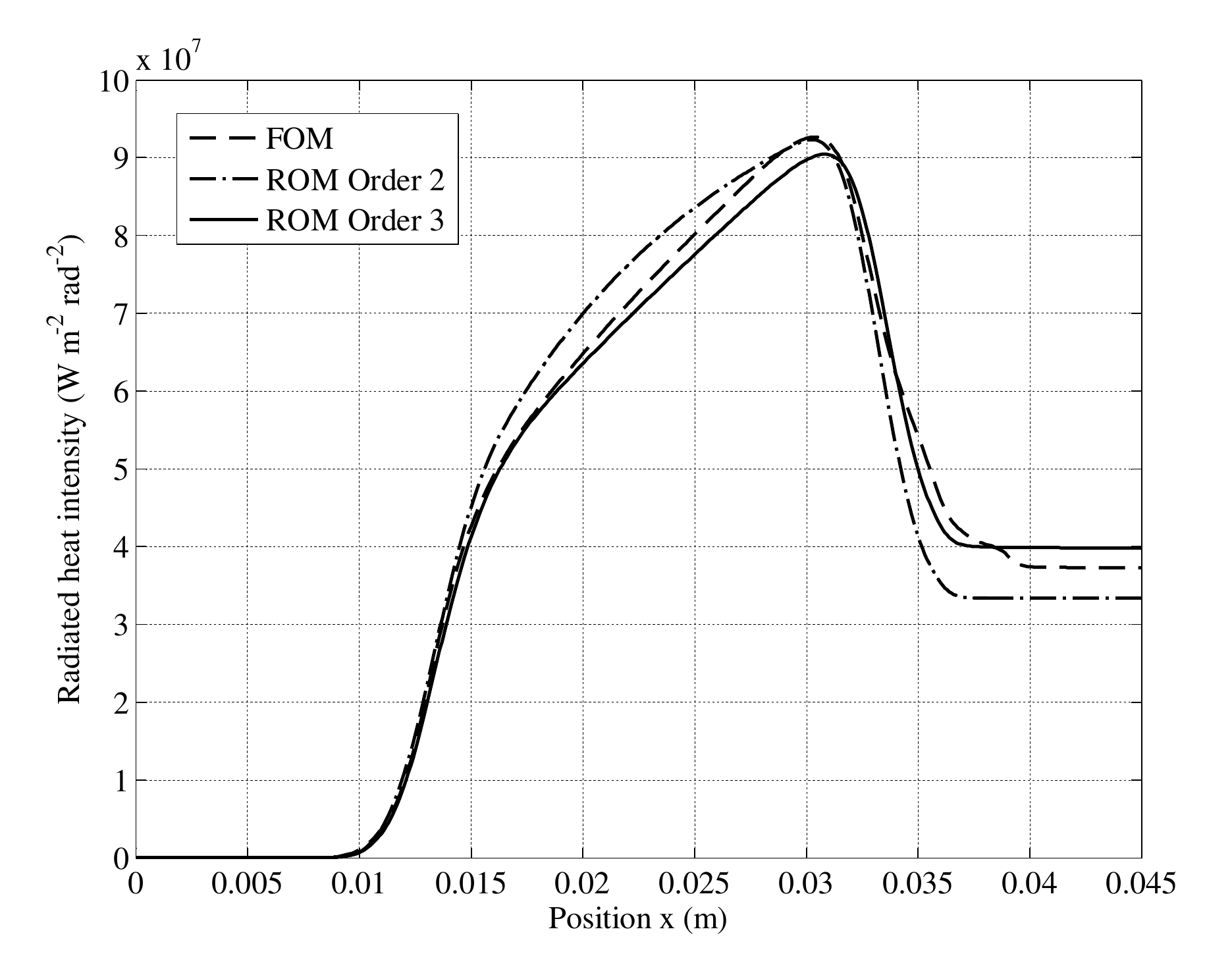}
\end{tabular}
} \protect\caption{\small (a) Temperature profile used to compare
the radiated heat predicted by the FOM with that predicted by the
ROM. (b) Comparison between the radiated heat intensities predicted
by the FOM and the ones given by two ROMs with order 2 and order 3,
respectively, using the temperature distribution of Figure (a).
Mixture: $25\%$ air, $50\%$ copper and $25\%$ hydrogen at
$2.5\,10^5\,$Pa.\normalsize}\label{F:FOM_ROM_compare}
\end{figure}
\section{Results}\label{S:results}

We applied our method to compute the equivalent absorption
coefficients and the partition of the EM spectrum in frequency bands
for several gaseous media, ranging from pure air to mixtures of
silver (or copper), air and hydrogen, and carbon dioxide and copper.
As an example, Figures \ref{F:abs_example1}-\ref{F:abs_example2}
show the results for pure air, using $M=2$ bands in the ROM. For
each temperature profile $\tilde{T}_j$, the weight $w_j$ in
\eqref{E:cost_term} has been chosen by taking the average value of
the corresponding intensity profile computed with the FOM. From
these results, it can be noted that for low pressure values (Figure
\ref{F:abs_example1}), in a temperature range of about
12-14$\,10^4\,$K the first of the two bands accounts for more than
80\% of the total black-body intensity, however with an absorption
coefficient of the order of $10^{-1}\,$m$^{-1}$. Considering that
the steady-state intensity is reached after about $3/\alpha\,$ of
propagation distance within the medium, this band can be considered
to be transparent with respect to the space-scale of interest (i.e.
about $10^{-2}\,$m). On the other hand, in the same temperature
range the second band has low emissivity (accounting for the
remaining 20\%), but large absorption coefficient, hence it is able
to reach in short distance the corresponding fraction of black-body
intensity. When large absorption coefficients are present (again
relative to the considered space-scale), the band is said to be
diffusive. In practice, however, the two bands emit (or absorb)
often similar values of intensity in the same distance, since the
transparent one contributes a small fraction of a very large
intensity value, while the diffusive one contributes a large
fraction of a small intensity value.

\begin{figure}[!hbt]
\centerline{
\begin{tabular}{cc}
(a)&(b)\\
\includegraphics[clip,width=6cm]{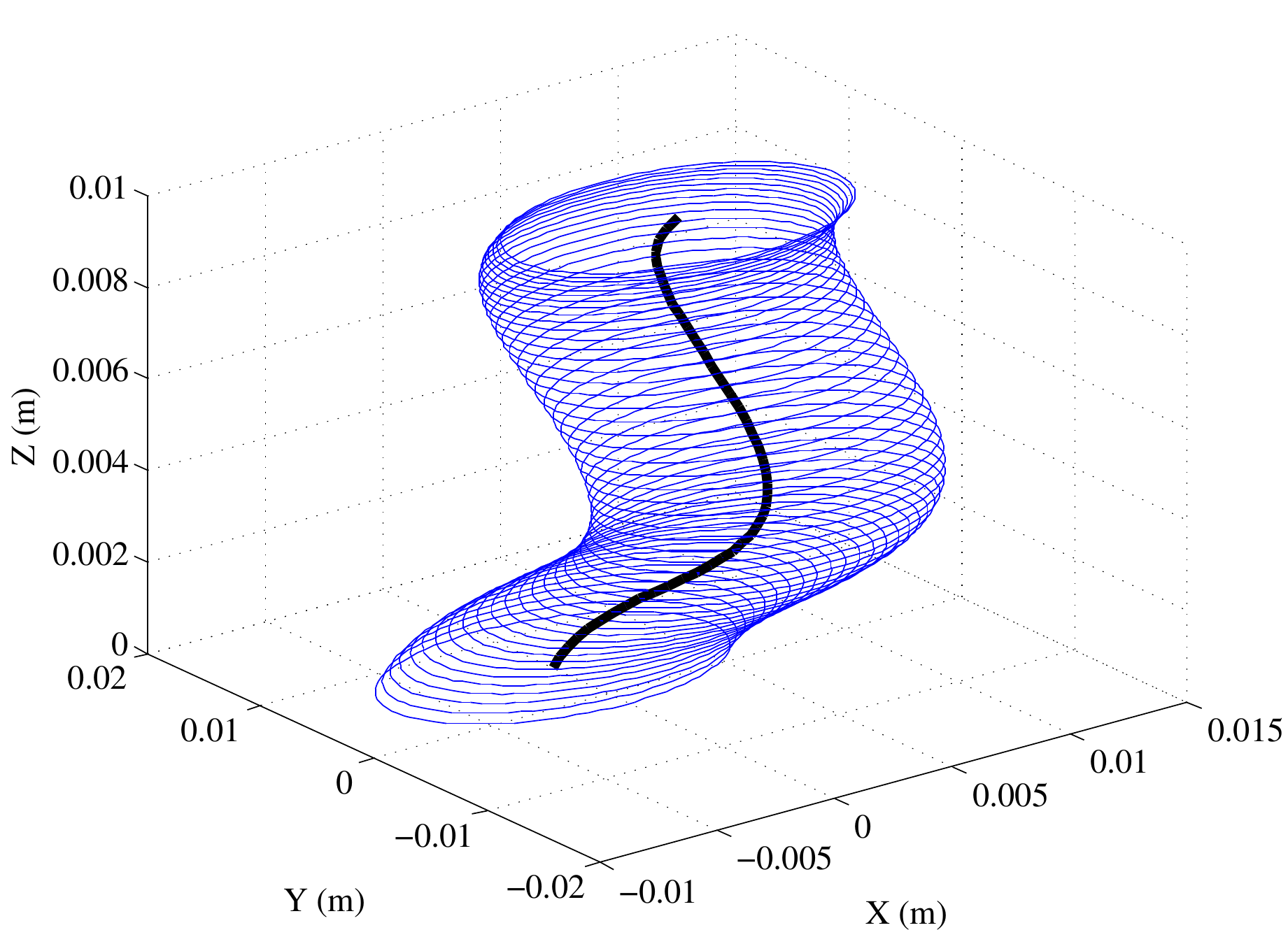}
&
\includegraphics[clip,width=6cm]{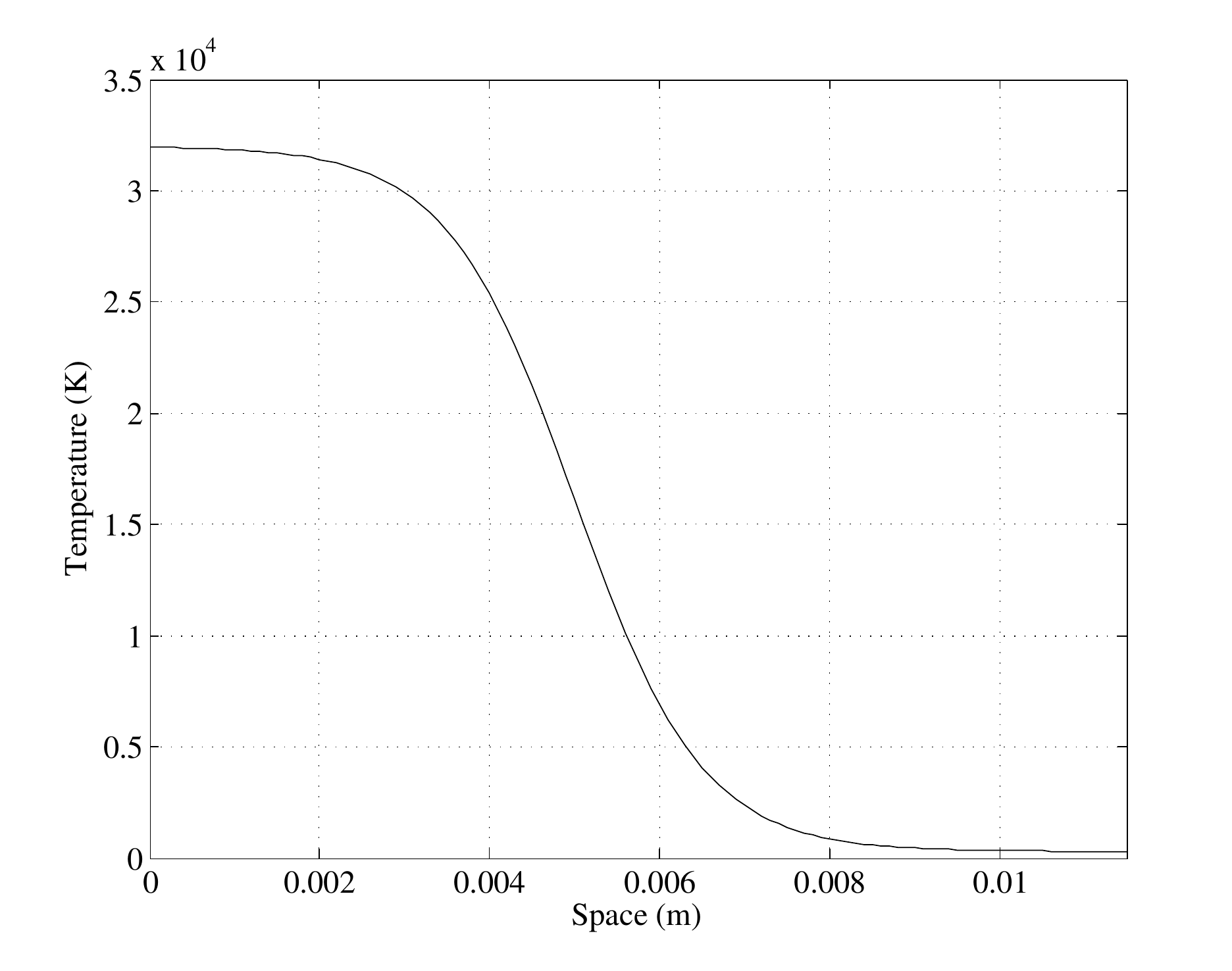}
\end{tabular}
} \protect\caption{\small Setup for a simulation of the radiated
heat intensity using the discrete ordinate method. A column of hot
gas ($75\%$ copper and $25\%$ air at $10^5\,$Pa) is standing in
front of a wall and we want to compute the radiated heat that
reaches the wall. (a) Scheme of the setup, the wall of interest is
parallel to the $(Y,Z)$ and contains the point
$(1.5\,10^{-2},0,0)\,$. The column is represented by blue circles
and the position of its center in $(X,Y,Z)$ is shown as a solid
black line. (b) Course of the temperature along a line parallel to
the $(X,Y)$ plane, passing through the center of the column, the
highest temperature value corresponding to the column
center.\normalsize}\label{F:sim_column}
\end{figure}
\begin{figure}[!hbt]
\centerline{
\includegraphics[clip,width=13cm]{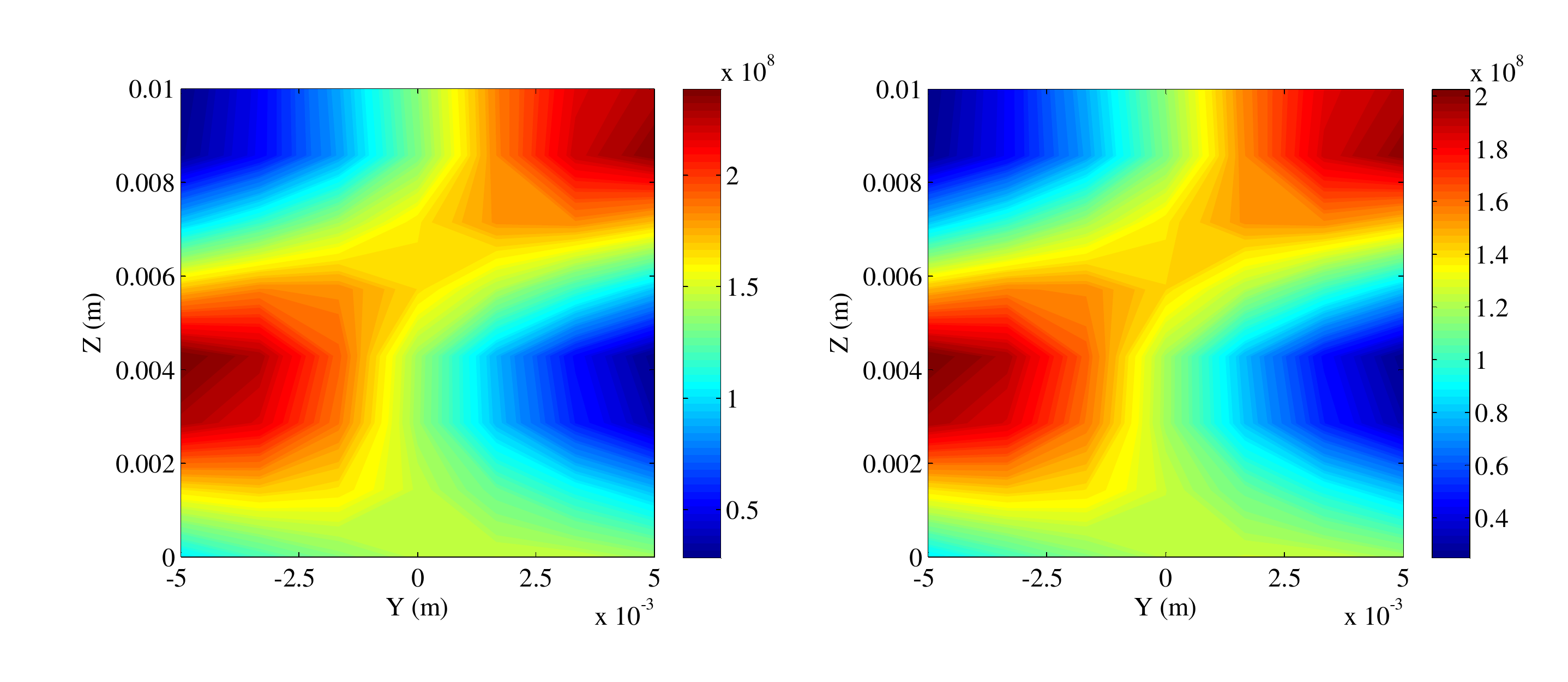}
} \protect\caption{\small Simulation results for the hot column
example with the DOM method. Radiated heat intensity (W/m$^2$) that
reaches the wall computed with the full-order model (left) and with
a reduced-order model (right) with three bands and piecewise affine
parametrization with 30 temperature nodes in the range
$[300,\,4\,10^4]\,$K.\normalsize}\label{F:sim_column_res}
\end{figure}

The situation can change significantly for higher pressure values:
in this case, the  diffusive bands can have larger emissivity values
at high temperature, meaning that they reach in short distance a
significant portion of the black-body intensity (see Figure
\ref{F:abs_example2}). This implies, roughly speaking, that at these
pressure and temperature values a larger quantity of radiated heat
is redistributed to the surroundings and transferred to the walls.
Similar considerations as the one just presented apply for all the
other mixtures that have been considered.

Figure \ref{F:Opt_bands}-(a) shows the dependency of the optimal
partition of the EM spectrum, i.e. of parameter
$\theta_{\hat{e}_{hf}}(\overline{p},\overline{y},\mu_1)$  on
pressure. It can be noted that, as pressure increases, the boundary
between the two bands shifts slightly and gradually towards lower
frequencies, meaning that, for a given temperature, the second band
accounts for larger fractions of the black-body intensity. This
effect is clear also from the courses of the emissivities in Figures
\ref{F:abs_example1}-\ref{F:abs_example2}. In Fig.
\ref{F:Opt_bands}-(b), the optimal partition achieved with $M=3$
bands is shown. It can be noted that the variation of the frequency
cuts with pressure is less marked in this case, in particular for
the boundary between the first and second bands. An example of the
corresponding absorption coefficients and emissivities is shown in
Fig. \ref{F:abs_example3}.

To give an idea of the accuracy achieved by the ROMs computed with
the proposed method with respect to the FOM, we present two further
examples. The first one is related to the radiated heat intensity
along a single line of propagation, with the temperature profile
shown in Figure \ref{F:FOM_ROM_compare}(a), through a mixture of
$25\%$ air, $50\%$ copper and $25\%$ hydrogen at $2.5\,10^5\,$Pa.
Such a profile is a realistic temperature distribution that can take
place in the plasma generated during the switching process of a
circuit breaker, and it was not part of the data used to derive the
ROM coefficients. Two ROMs are considered, one with two bands and
the other with three bands. The resulting distributions of the
radiated heat intensity are shown in Figure
\ref{F:FOM_ROM_compare}(b). It can be noted that both ROMs are able
to reproduce quite accurately the intensity profile given by the
FOM, with the third-order one being slightly more accurate. This
example thus shows that increasing the number of bands gives in
general higher accuracy, but the gain is smaller and smaller (in
line with the considerations of section \ref{S:FreqDom}) and usually
it is not worth using more than three bands, due to the increased
computational load in the CFD simulations.

The second example is related to the use of the derived ROM for a
discrete-ordinate method (DOM, see e.g. \cite{SiHo92}) simulation,
where we want to compute the radiated heat intensity that reaches a
wall in front of a column of hot gas composed by $75\%$ copper and
$25\%$ air, at $10^5\,$Pa. Such a setup is described in Figure
\ref{F:sim_column}(a)-(b).

This example is more meaningful for the sake of CFD simulations of a
plasma in real devices, where not just the intensity along a single
line but the net total intensity obtained by integrating over all
directions has to be computed, for each one of the finite volumes or
elements that partition the computational domain. The intensity
distribution (in W/m$^2$) on the wall is shown in Figure
\ref{F:sim_column_res}. It can be noted that the ROM is able to
capture well the qualitative behavior given by the FOM, with an
average error of about $10\%$. Such a value is indeed very small as
compared with the accuracy that can be achieved with other
approaches, e.g. using constant absorption coefficients over
frequency bands.

Finally, as regards the computational load, here are some
indications obtained by running the order reduction algorithm on a
workstation equipped with 12 Intel Xeon$^\circledR$ cores at
2.4$\,$GHz each and 52GB of RAM in total. The average time required
to compute one set of coefficients for a given pressure and
composition pair $(\overline{p},\overline{y})$ was  0.23 hours for
order two models, and 0.32 hours for order three models, in both
cases with a temperature discretization with 50 nodes for the
absorption coefficient, i.e. $101$ and $152$ optimization variables,
respectively (see section \ref{SSS:parametrization}).

\section{Conclusions and future developments}\label{S:conclusions}
We presented a new approach to derive reduced order models of the
radiated heat intensity through a participating gaseous medium. The
approach is based on nonlinear identification and numerical
optimization methods, and it can deliver models with low complexity
but still highly accurate with respect to the original full
absorption spectrum of the considered mixture. We also introduced a
system-theoretic perspective of the phenomenon, as well as a
frequency-domain analysis that justifies the use of low-order models
to predict the total radiated intensity. We presented several
comparisons between the full-order model and the reduced-order ones,
obtained with our approach, and discussed the required computational
effort. With the proposed technique, one can reliably and
systematically derive radiation models that can be used in CFD
simulations of a plasma, since both the average absorption
coefficients and the frequency bands are a result of the numerical
optimization. The accuracy of the model depends mainly on the
accuracy of the base data that the model approximates. Future
development efforts shall thus be devoted to assess and improve such
base data.

\section*{Acknowledgements}
The authors would like to thank Prof. V. Aubrecht for providing the
the absorption spectra used as base data in this study.

\section*{References}
\bibliographystyle{plain}
\bibliography{Biblio}




\end{document}